%Dehn surgery and negatively curved 3-manifolds

\input epsf
\input amssym.def
\input amssym
\magnification=1100
\baselineskip = 0.25truein
\lineskiplimit = 0.01truein
\lineskip = 0.01truein
\vsize = 8.5truein
\voffset = 0.7truein
\parskip = 0.10truein
\parindent = 0.3truein
\settabs 10 \columns
\hsize = 5.4truein
\font\ninerm=cmr9
\font\nineit=cmti9

\hoffset=0.5truein

\def\sqr#1#2{{\vcenter{\vbox{\hrule height.#2pt
\hbox{\vrule width.#2pt height#1pt \kern#1pt
\vrule width.#2pt}

\hrule height.#2pt}}}}
\def\square{\mathchoice\sqr46\sqr46\sqr{3.1}6\sqr{2.3}4}

\font\bigtenrm=cmr10 scaled 1400

\vskip 2in
\centerline{\bf {\bigtenrm DEHN SURGERY AND}}
\centerline{\bf {\bigtenrm NEGATIVELY CURVED 3-MANIFOLDS}}
\tenrm
\vskip 14pt
\centerline{DARYL COOPER AND MARC LACKENBY}
\vskip 10pt

\tenrm

\footnote{}{\ninerm 1991 {\nineit Mathematics Subject Classification.}
Primary 57N10; Secondary 57M25.}

\vskip 18pt
\centerline {1. INTRODUCTION}
\vskip 6pt

Dehn surgery is perhaps the most common way of constructing 3-manifolds,
and yet there remain some profound mysteries about its behaviour.
For example, it is still not known whether there exists a
3-manifold which can be obtained from $S^3$
by surgery along an infinite number of distinct knots.
\baselineskip=0.10truein \footnote{$^1$}{\ninerm
Since this paper was written, John Osoinach has constructed a family of
3-manifolds, each with infinitely many knot surgery descriptions
[Ph.D. Thesis, University of Texas at Austin].}
\baselineskip=0.25truein
(See Problem 3.6 (D) of Kirby's list [9]). In this paper,
we offer a partial solution to this problem, and exhibit
many new results about Dehn surgery. The methods we employ make use of
well-known constructions of negatively curved metrics on certain
3-manifolds. 

We use the following standard terminology. A {\sl slope}
on a torus is the isotopy class of an unoriented essential simple closed curve.
If $s$ is a slope on a torus boundary component of a 3-manifold $X$,
then $X(s)$ is defined to be the 3-manifold
obtained by Dehn filling along $s$. More generally,
if $s_1, \dots, s_n$ is a collection of slopes on distinct
toral components of $\partial X$, then we write $X(s_1, \dots, s_n)$
for the manifold obtained by Dehn filling along each of these slopes.

We also abuse terminology in the standard way
by saying that a compact orientable
3-manifold $X$, with $\partial X$ a (possibly empty)
union of tori, is {\sl hyperbolic}
if its interior has a complete finite volume hyperbolic structure.
If $X$ is hyperbolic, we also say that the core of the filled-in solid
torus in $X(s)$ is a {\sl hyperbolic knot}.

Let $X$ be a hyperbolic 3-manifold,
and let $T_1, \dots, T_n$ be a collection of components
of $\partial X$. Now, associated with each torus $T_i$,
there is a cusp in ${\rm int}(X)$ homeomorphic to $T^2 \times [1, \infty)$.
We may arrange that the $n$ cusps are all disjoint.
They lift to an infinite set of disjoint
horoballs in ${\Bbb H}^3$. Expand these horoballs equivariantly
until each horoball just touches some other. Then, the
image under the projection map ${\Bbb H}^3 \rightarrow {\rm int}(X)$
of these horoballs is a {\sl
maximal horoball neighbourhood} of the cusps at $T_1, \dots, T_n$.
When $n =1$, this maximal horoball neighbourhood is unique.
Let ${\Bbb R}_i^2$ be the boundary in ${\Bbb H}^3$ of one of these
horoballs associated
with $T_i$. Then ${\Bbb R}_i^2$
inherits a Euclidean metric from ${\Bbb H}^3$.
A slope $s_i$ on $T_i$ determines a primitive element
$[s_i] \in \pi_1(T_i)$, which is defined up to sign. This
corresponds to a covering translation of ${\Bbb R}_i^2$,
which is just a Euclidean translation.
We say that $s_i$ has
{\sl length} $l(s_i)$ given by the length of the
associated translation vector.
When $n > 1$, this may depend upon the choice of a maximal
horoball neighbourhood of $T_1 \cup \dots \cup T_n$, but
for $n = 1$, the length of $s_1$ is a topological
invariant of the manifold $X$ and the slope $s_1$,
by Mostow Rigidity [Theorem C.5.4, 3].
The concept of slope length is very relevant to Dehn surgery along
hyperbolic knots, and plays a crucial r\^ole in this paper.
Note that slope length is measured in the metric on $X$, not in
any metric that $X(s_1, \dots, s_n)$ may happen to have.
This notion of slope length arises in the following well-known theorem of
Gromov and Thurston, the so-called `$2 \pi$' theorem.

\noindent {\bf Theorem} (Gromov, Thurston [Theorem 9, 4]).
{\sl Let $X$ be a compact
orientable hyperbolic 3-manifold. Let $s_1, \dots, s_n$ be a
collection of slopes on distinct components $T_1, \dots, T_n$ of
$\partial X$. Suppose that there is a horoball neighbourhood of
$T_1 \cup \dots \cup T_n$ on which each $s_i$ has length greater than $2 \pi$.
Then $X(s_1, \dots, s_n)$ has a
complete finite volume Riemannian metric with all sectional curvatures
negative.}

The following theorem, which is the main result of this paper,
asserts that (roughly speaking) any given 3-manifold $M$ can be constructed
in this fashion in at most a finite number of ways.

\noindent {\bf Theorem 4.1.} {\sl Let $M$ be a compact
orientable 3-manifold, with
$\partial M$ a (possibly empty) union of tori.
Let $X$ be a hyperbolic manifold and let $s_1, \dots, s_n$
be a collection of slopes on $n$ distinct tori $T_1, \dots, T_n$
in $\partial X$,
such that $X(s_1, \dots, s_n)$ is homeomorphic to $M$.
Suppose that there exists in ${\rm int}(X)$ a maximal horoball neighbourhood
of $T_1 \cup \dots \cup T_n$
on which each slope $s_i$ has length at least
$2 \pi + \epsilon$, for some $\epsilon >0$.
Then, for any given $M$ and $\epsilon$, there is
only a finite number of possibilities (up to isometry)
for $X$, $n$ and $s_1, \dots, s_n$.}

This is significant because `almost all' closed orientable 3-manifolds are
obtained by such a Dehn surgery.
More precisely, any closed orientable 3-manifold is
obtained by Dehn filling some hyperbolic 3-manifold $X$ ([13] and [12]).
After excluding at most $48$
slopes from each component of $\partial X$, all remaining
slopes have length more than $2 \pi$ [Theorem 11, 4].

Theorem 4.1 has the following corollary,
which is a partial solution to Kirby's Problem 3.6 (D).

\noindent {\bf Corollary 4.5.} {\sl For a given closed orientable
3-manifold $M$, there is at most a finite number of hyperbolic knots $K$ in $S^3$
and fractions $p/q$ (in their lowest terms) such that $M$ is
obtained by $p/q$-Dehn surgery along $K$ and $\vert q \vert > 22$.}

Dehn filling also arises naturally in the study of branched
covers. Recall [15] that a branched cover of a 3-manifold $Y$ over
a link $L$ is determined by a transitive representation
$\rho \colon \pi_1(Y - L) \rightarrow S_r$, where $S_r$ is the
symmetric group on $r$ elements. The stabiliser of one of these
elements is a subgroup of $\pi_1(Y-L)$
which determines a cover $X$ of $Y - {\rm int}({\cal N}(L))$.
The branched cover is then obtained by Dehn filling
each component $P$ of $\partial X$ that is a lift of
some component of $\partial {\cal N}(L)$. The Dehn filling
slope on $P$ is the slope which the covering map sends
to a multiple of a meridian
slope on $\partial {\cal N}(L)$. This multiple is known
as the {\sl branching index} of $P$.
Branched covers are a surprisingly general construction. For example,
any closed orientable 3-manifold is a branched cover of
$S^3$ over the figure-of-eight knot [7]. Thus, the
following corollary to Theorem 4.1 is useful.

\noindent {\bf Corollary 4.8.} {\sl Let $M$ be a compact
orientable 3-manifold with $\partial M$ a (possibly empty)
union of tori, which is obtained as a branched
cover of a compact orientable 3-manifold $Y$ over a hyperbolic link $L$,
via representation $\rho \colon \pi_1(Y-L) \rightarrow S_r$.
Suppose that the branching index of every lift
of every component of $\partial {\cal N}(L)$ is
at least 7. Then, for a given $M$, there are only finitely
many possibilities for $Y$, $L$, $r$ and $\rho$.}

This paper is organised as follows. In Section 2, we establish lower
bounds on slope length from topological information. In Section 3,
we review the proof of the `$2 \pi$' theorem
and establish a `controlled' version of the theorem, which involves
estimates of volume and curvature. In Section 4, the main results
about Dehn surgery are deduced from the work in Sections 2 and 3.
In Section 5, we obtain restrictions on the genus of
surfaces in the complement of a hyperbolic knot in terms of their
boundary slopes. In Section 6, we examine 3-manifolds which
are `almost hyperbolic', in the sense that for any $\delta > 0$,
they have a complete finite volume 
Riemannian metric with all sectional curvatures
between $-1 - \delta$ and $-1+\delta$. We show that such 3-manifolds
must have a complete finite volume hyperbolic structure.
The proof of this result uses Dehn surgery in a crucial way.

\vskip 18pt
\centerline {2. WHEN IS A SLOPE LONG?}
\vskip 6pt

Since the majority of theorems in this paper are stated in terms of
slope length, we will now establish some conditions which imply
that a slope is long. Throughout this section, we will
examine slopes lying on a single torus $T$ in $\partial X$.
There may be boundary components of $X$ other than $T$,
but nevertheless the length of a slope on $T$ is a 
well-defined topological invariant.

Recall that the {\sl distance} $\Delta(s_1, s_2)$
between two slopes $s_1$ and $s_2$ on a torus
is defined to be the minimum number of intersection points
of two representative simple closed curves.
The following lemma implies that if the distance between two slopes
is large, then at least one of them must be long.

\noindent {\bf Lemma 2.1.} {\sl Let $X$ be a hyperbolic 3-manifold
with a torus $T$ in its boundary. Let $s_1$ and $s_2$ be slopes
on $T$. Then
$$l(s_1) \, l(s_2) \geq \sqrt 3 \, \Delta(s_1, s_2).$$
Moreover, if all slopes on $T$ have length at least $L$, say, then
$$l(s_1) \, l(s_2) \geq \sqrt 3 \, L^2 \Delta(s_1, s_2).$$}

\noindent {\sl Proof.} 
It is a well-known observation of Thurston that the length
of each slope on $T$ is at least 1. (See [Theorem 11, 4] for example.) Hence,
the first inequality of the lemma follows from the
second. We pick a generating set for $H_1(T)$ as follows.
We first assign arbitrary orientations to $s_1$ and $s_2$.
We let $[s_1] \in H_1(T)$ be one generator, and extend this to
a generating set by picking one further element $[s_3]$.
Then,
$$[s_2] = \pm \Delta(s_1,s_2) [s_3] + n [s_1],$$
for some integer $n$.
Let $N$ be the maximal horoball neighbourhood of the cusp at $T$.
Let ${\Bbb R}^2$ be the boundary in ${\Bbb H}^3$ of an associated
horoball. Let $P$ (respectively, $P'$)
be a fundamental domain in ${\Bbb R}^2$ for the
group of covering translations generated by $[s_1]$ and $[s_2]$
(respectively, by $[s_1]$ and $[s_3]$). Note that
$$l(s_1) \, l(s_2) \geq {\rm Area}(P) =
\Delta(s_1, s_2) {\rm Area}(P').$$
This formula is clear from Fig. 1.

\vskip 24pt
\epsfxsize = 3in
\centerline{\epsfbox{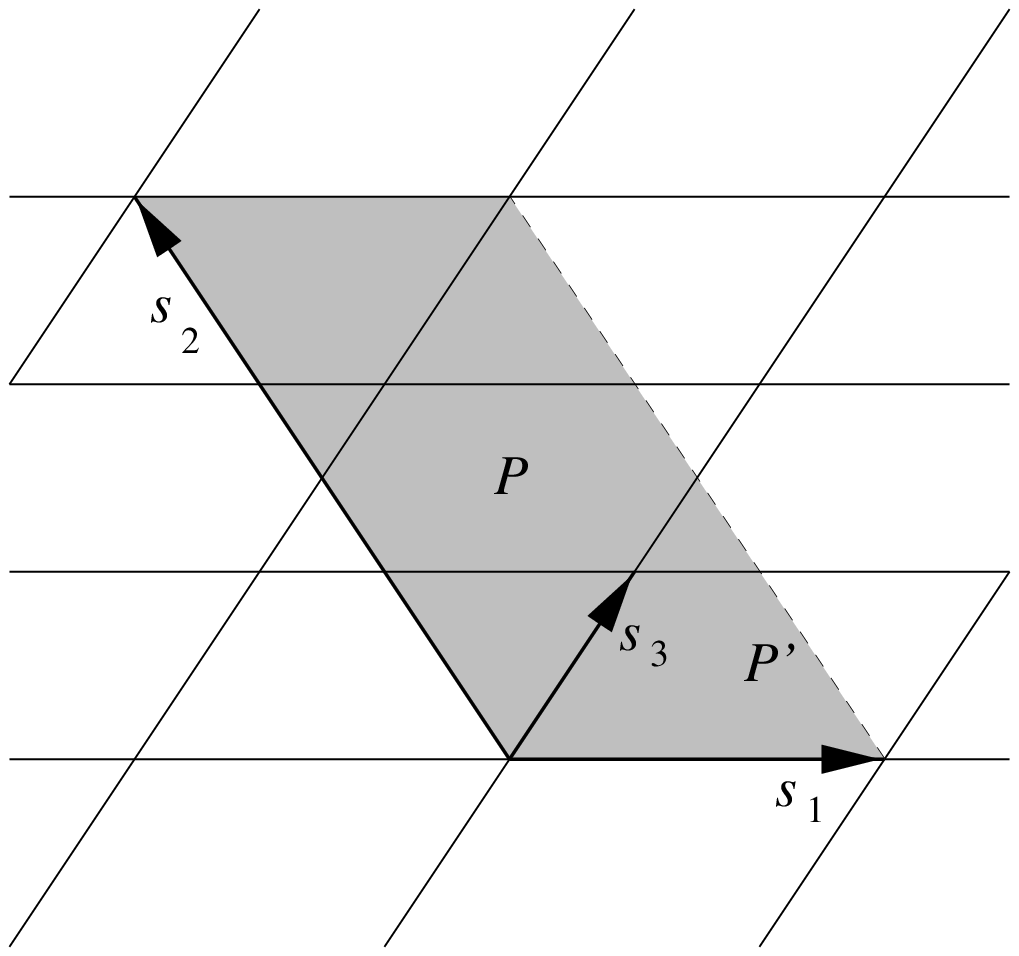}}
\centerline{Fig. 1}
\vskip 12pt

It is well known that the area of $P'$ is at least $\sqrt 3$ (see
[Theorem 2, 1]). However, the argument there readily implies that
$${\rm Area}(P') \geq \sqrt 3 \, L^2.$$
Hence, we deduce that
$$l(s_1) \, l(s_2) \geq \sqrt 3 \, L^2 \Delta(s_1, s_2).\eqno \square$$

Now, there are various topological circumstances when a slope $e$
is known to be `short'. These are summarised in the following
proposition.

\noindent {\bf Proposition 2.2.} {\sl Let $X$ be a compact
orientable hyperbolic 3-manifold and
let $T$ be a toral boundary component of $X$. Then a slope $e$ on $T$
has length no more than $2 \pi$ if either of the following hold:

\item{1.} ${\rm int}(X(e))$ does not admit a complete finite volume
negatively curved Riemannian metric (for example,
$X(e)$ may be reducible, toroidal or Seifert fibred), or

\item{2.} the core of the filled-in solid torus in $X(e)$ has finite
order in $\pi_1(X(e))$.}

\noindent {\sl Proof.} Part (1) above is a mere restatement of the
`$2 \pi$' theorem, with the added assertion that if the
interior of a compact orientable
3-manifold $M$ admits a complete finite volume negatively curved
metric, then $M$ cannot be reducible, toroidal or Seifert fibred.
This is well-known, but we sketch a proof. By the Hadamard-Cartan
theorem [2], the universal cover of ${\rm int}(M)$ is homeomorphic to ${\Bbb R}^3$,
and so $M$ is irreducible.  By [2],
any ${\Bbb Z} \oplus {\Bbb Z}$ subgroup of $\pi_1(M)$
is peripheral. Hence, $M$ is atoroidal.
Since ${\rm int}(M)$ is covered by ${\Bbb R}^3$, $\pi_1(M)$ is infinite.
The fundamental group of
any orientable Seifert fibre space either has non-trivial centre
or is trivial [16]. However, the fundamental group of a complete
negatively curved finite volume Riemannian manifold cannot
have non-trivial centre [2]. Hence, $M$ is not Seifert fibred.

To prove Part (2), we recall from the
proof of the `$2 \pi$' theorem [4] that if $l(e) > 2\pi$, then
${\rm int}(X(e))$ has a complete finite volume negatively curved metric,
in which the core of the filled-in solid torus is a geodesic.
But, in such a manifold,
closed geodesics have infinite order in the fundamental group [2].
$\square$

We therefore make the following definition.

\noindent {\bf Definition 2.3.} Let $X$ be a hyperbolic 3-manifold and
let $e$ be a slope on a toral boundary component $T$ of $X$.
Then $e$ is {\sl short} if $l(e) \leq 2 \pi$. We say that $e$ is
{\sl minimal} if $l(s) \geq l(e)$ for all slopes $s$ on $T$.

Note that there is always at least one minimal slope on $T$, but
that it need not be short.
The importance of the above definition is that if some slope $s$
has large intersection number with a slope $e$ which is either short
or minimal, then we can deduce that the length of $s$ is large. Moreover
the bounds we construct are independent of the manifold $X$.

\noindent {\bf Corollary 2.4.} {\sl Let $X$ be a compact hyperbolic 3-manifold,
and let $s$ be a slope on a torus component $T$ of
$\partial X$. If $e$ is a short slope on $T$, then
$$l(s) \geq \sqrt 3 \, \Delta(s,e) / 2\pi.$$
If $e$ is a minimal slope on $T$, then
$$l(s) \geq \sqrt 3 \, \Delta(s,e).$$
}

\noindent {\sl Proof.} The first inequality is an immediate consequence
of Lemma 2.1 and the definition of `short'. If $e$ is a minimal slope on $T$,
then all slopes on $T$ have length
at least $l(e)$, and hence by Lemma 2.1,
$$l(s) \, l(e) \geq \sqrt 3 \, [l(e)]^2 \Delta(s, e).$$
Thus, we get that
$$l(s) \geq \sqrt 3 \, l(e) \Delta (s,e) \geq \sqrt 3 \, \Delta(s,e),$$
since $l(e) \geq 1$ by [Theorem 11, 4]. $\square$

\vskip 18pt
\centerline {3. ESTIMATES OF CURVATURE, VOLUME AND GROMOV NORM}
\vskip 6pt

In this section, we compare the Gromov norm [Chapter 6, 17]
of a hyperbolic 3-manifold $X$
with the Gromov norm of a 3-manifold obtained by
Dehn filling tori $T_1, \dots, T_n$ in $\partial X$. 

In [Section 6.5, 17], Thurston gave various definitions of
the Gromov norm of a compact orientable
3-manifold $X$, with $\partial X$ a (possibly empty) union of tori.
We shall use the terminology $\vert X \vert$ for the
quantity which Thurston calls $\vert \vert [X, \partial X] \vert \vert_0$.
This is defined as follows. Consider the fundamental class 
$[X, \partial X]$ in the singular homology group 
$H_3(X, \partial X; {\Bbb R})$. If $z = \sum a_i \sigma_i$ is
a representative of $[X, \partial X]$, where $a_i \in {\Bbb R}$
and each $\sigma_i$ is a singular 3-simplex, then we
consider the real number $\vert \vert z \vert \vert = \sum \vert a_i \vert$.
In the case where $\partial X = \emptyset$, the Gromov
norm is defined to be
$$\vert X \vert = \inf \lbrace \vert \vert z \vert \vert : z 
\ \hbox{represents} \ [X, \partial X] \rbrace.$$
In the case where $\partial X \not= \emptyset$, a
representative $z$ of $[X, \partial X]$ determines
a representative $\partial z$ of $[\partial X] \in 
H_2(\partial X; {\Bbb R})$. Thurston defines
$$\vert X \vert = \liminf_{a \rightarrow 0} \lbrace
\vert \vert z \vert \vert : z \ \hbox{represents}
\ [X, \partial X] \ \hbox{and} \ \vert \vert \partial z \vert \vert \leq
\vert a \vert \rbrace$$
and shows that this limit exists. Crucially, the
Gromov norm of $X$ is a topological invariant.

We compare the Gromov norms of $X$ and $X(s_1, \dots, s_n)$ (where
$s_1, \dots, s_n$ are slopes on $\partial X$) by
returning to the proof of the `$2 \pi$' theorem.
The idea behind this proof is simple.
One first removes from $X$ the interior of an almost maximal
horoball neighbourhood $N$ of the cusps at $T_1 \cup \dots \cup T_n$.
Then one glues back in
solid tori $V_i$ which have negatively curved
Riemannian metrics agreeing near
$\partial V_i$ with that near $\partial N$. The
following proposition deals with the sectional curvatures
and the volume of the metric on each solid
torus.

If $M$ is a manifold with interior having a Riemannian metric $k$,
let ${\rm Vol}(M,k)$
denote its volume and let $\kappa_{\rm inf}(M,k)$ (respectively,
$\kappa_{\rm sup}(M,k)$) denote the infimum (respecively, the supremum) of
its sectional curvatures.

\noindent {\bf Proposition 3.1.} {\sl For any two real numbers $\ell_1 > 2\pi$
and $\ell_2 > 0$, we may construct a Riemannian
metric $k$ on the solid torus $V$, with the following properties.
In a collar neighbourhood of $\partial V$, the metric is hyperbolic.
The boundary $\partial V$ inherits a Euclidean metric
$k\vert_{\partial V}$. The length in this metric
of a shortest meridian curve $C$ on $\partial V$
is $\ell_1$. The length of a (Euclidean) geodesic running
perpendicularly from $C$ to $C$ is $\ell_2$.
Also, ${\rm Vol}(V,k) / {\rm Vol}(\partial V, k \vert_{\partial V})$,
$\kappa_{\rm inf}(V,k)$ and $\kappa_{\rm sup}(V, k)$ are all
independent of $\ell_2$.
But, there is a non-decreasing function $\alpha \colon (2 \pi , \infty)
\rightarrow (0,1)$ such that
$$-(\alpha(\ell_1))^{-1} \leq \kappa_{\rm inf}(V,k) <
\kappa_{\rm sup}(V,k) \leq -\alpha(\ell_1),$$
$${{\rm Vol}(V,k) \over
{\rm Vol}(\partial V,k \vert_{\partial V})/2} \geq \alpha(\ell_1).$$
}

\noindent {\sl Proof.} In Bleiler and Hodgson's proof of the
`$2 \pi$' theorem [4], a Riemannian metric $k$ is constructed on $V$
which has most of these properties.
They assign cylindrical co-ordinates $(r, \ \mu, \ \lambda)$
to $V$, where $r \leq 0$ is the radial distance
measured {\it outward} from  $\partial V$, $0 \leq \mu \leq 1$ is
measured in the meridional direction and $0 \leq \lambda \leq 1$
is measured in a direction perpendicular to $\mu$ and $r$.
The distance from the core of $V$ to the boundary is $-r_0$,
for some negative constant $r_0$. The Riemannian metric
is denoted
$$ds^2 = dr^2 + [f(r)]^2 d\mu^2 + [g(r)]^2 d \lambda^2,$$
where $f \colon [r_0,0] \rightarrow {\Bbb R}$
and $g \colon [r_0,0] \rightarrow {\Bbb R}$ are functions.
Graphs of $f$ and $g$ are given in Bleiler and Hodgson's paper [4],
and are reproduced in Fig. 2.

\vskip 24pt
\epsfxsize = 4.5in
\centerline{\epsfbox{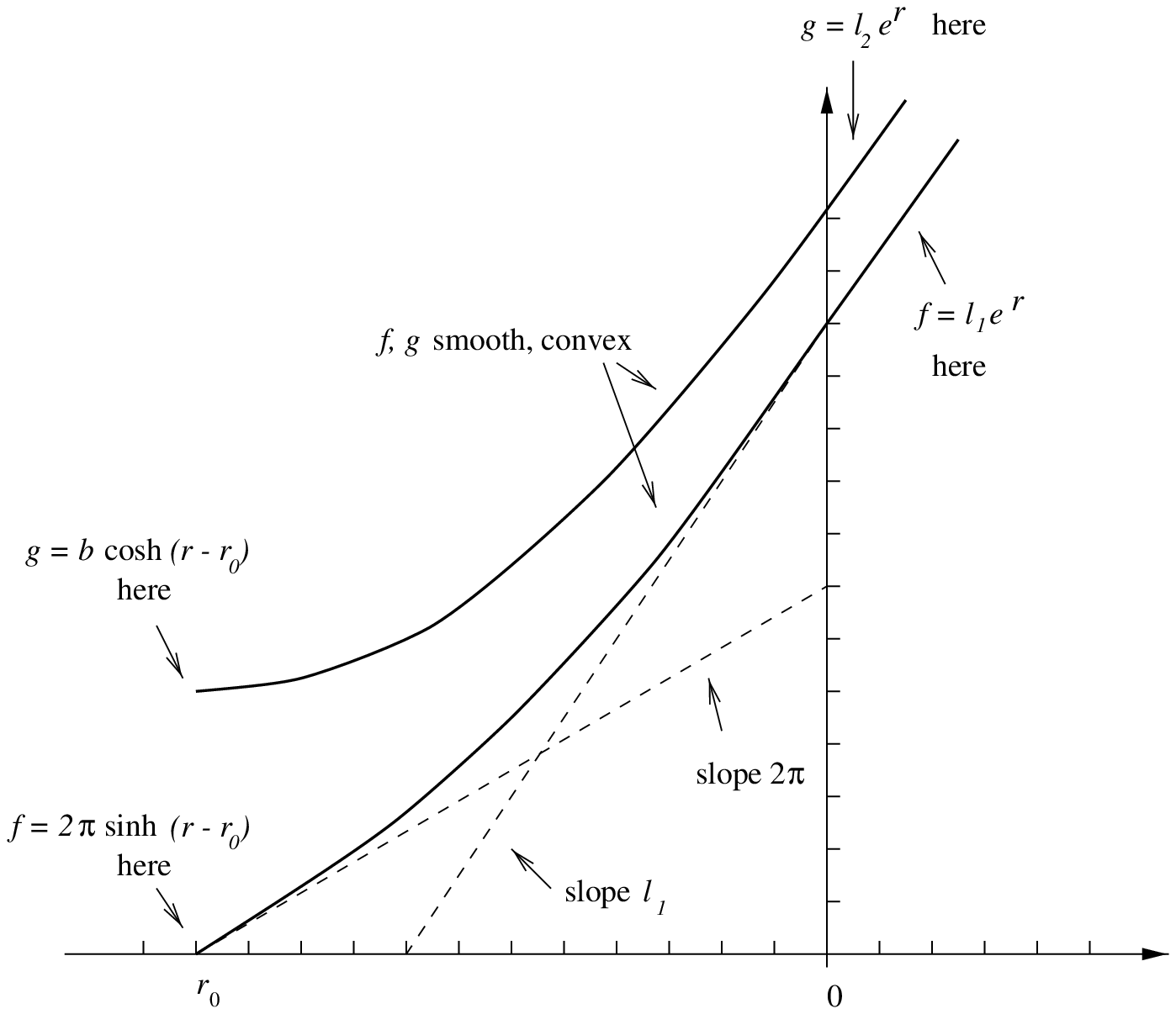}}
\centerline{Fig. 2}
\vskip 12pt

Making the substitutions $x_1 = r$, $x_2 = \mu$ and
$x_3 = \lambda$, they calculate the sectional curvatures as
$$\eqalign{\kappa_{12} &= -{f'' \over f},\cr
\kappa_{13} &= -{g'' \over g},\cr
\kappa_{23} &= -{f' \cdot g' \over f \cdot g}.\cr}$$
They observe that, since $f$, $f'$, $f''$, $g$, $g'$
and $g''$ are all positive in the range $r_0 < r \leq 0$, then the
sectional curvatures are all negative. 
To ensure that the cone angle at the core of $V$ is
$2 \pi$, it is enough to ensure that the gradient of
$f$ at $r = r_0$ is $2 \pi$. Also, near $r = 0$,
$f$ and $g$ are both exponential, which guarantees
that the sectional curvatures are all $-1$ near
$\partial V$.

Bleiler and Hodgson argue that, providing $\ell_1 > 2 \pi$
and $\ell_2 > 0$, we may find a metric $k$ satisfying
the above properties. We can therefore pick $\alpha$ (which is a 
real-valued function
of $\ell_1$ and $\ell_2$) as follows. Given $\ell_1$ and
$\ell_2$, there is a Riemannian metric $k$ on $V$ satisfying
all of the above properties, and a real number $a$
(with $0 < a < 1$) for which the following inequalities
hold:
$$-a^{-1} < \kappa_{\rm inf}(V, k) < \kappa_{\rm sup}(V, k)
< -a,$$
$${{\rm Vol}(V,k) \over
{\rm Vol}(\partial V,k \vert_{\partial V})/2} > a.$$
We define $\alpha(\ell_1, \ell_2)$ to be {\sl half}
the supremum of $a$, where $a$ satisfies the above
inequalities and $k$ satisfies the above conditions.
There is a good deal of freedom over the choice
of $\alpha$. We picked half the supremum simply
because it is less than the supremum, and hence
there is some metric $k$ for which
$$-(\alpha(\ell_1, \ell_2))^{-1} < \kappa_{\rm inf}(V, k) 
< \kappa_{\rm sup}(V, k) < -\alpha(\ell_1, \ell_2),$$
$${{\rm Vol}(V,k) \over
{\rm Vol}(\partial V,k \vert_{\partial V})/2} > 
\alpha(\ell_1, \ell_2).$$

At this stage, $\alpha$ depends on both $\ell_1$
and $\ell_2$. But we shall show that $\alpha$ is
independent of $\ell_2$ and is a non-decreasing function
of $\ell_1$. 

Given a metric $k$ on $V$
$$ds^2 = dr^2 + [f(r)]^2 d\mu^2 + [g(r)]^2 d \lambda^2,$$
we can define another metric
$$ds_0^2 = dr^2 + [f(r)]^2 d\mu^2 + c^2 [g(r)]^2 d\lambda^2,$$
for any positive real constant $c$. Using the formulae for
the sectional curvatures, we see that this alteration leaves the
sectional curvatures unchanged. It does not alter
the ratio of ${\rm Vol}(V, k)$
and ${\rm Vol}(\partial V, k \vert_{\partial V})$.
It leaves $\ell_1$ unchanged, but scales
$\ell_2$ by a factor of $c$. Hence,
$\alpha$ is independent of $\ell_2$, and
we therefore refer to $\alpha$ as a function of the
single variable $\ell_1$.
Note that we could not have made a similar argument
with $\ell_1$, since it is vital that the cone angle
at the core of $V$ is $2 \pi$.

It remains to show that $\alpha$ is a non-decreasing function
of $\ell_1$. Suppose $V$ has a metric $k$ as above. Then
we can enlarge $V$ (creating a bigger solid
torus $V'$ with metric $k'$) by letting $r$ vary in the
range $r_0 \leq r \leq c$, for
any positive real constant $c$, and defining
$f(r) = \ell_1 e^r$ and $g(r) = \ell_2 e^r$
for $r \geq 0$. In other words, we attach a collar
to $\partial V$, with the metric being hyperbolic in
the collar. The length of the shortest meridian
curve on the boundary of the new solid torus $V'$ is $\ell_1 e^c$.
If the metric $k$ on $V$ satisfies
$$-a^{-1} < \kappa_{\rm inf}(V, k) < \kappa_{\rm sup}(V, k)
< -a,$$
$${{\rm Vol}(V,k) \over
{\rm Vol}(\partial V,k \vert_{\partial V})/2} > a,$$
(with $0 < a < 1$), then the metric on the enlarged solid torus satisfies the
same inequalities. The inequality regarding volumes
requires some explanation:
$$\eqalign{{\rm Vol}(V',k') &= {\rm Vol}(V, k)
+ \int_0^c \, f(r) g(r) dr \cr
&= {\rm Vol}(V,k) + \ell_1 \ell_2 (e^{2c} - 1) /2 \cr
&> a {\rm Vol}(\partial V,k \vert_{\partial V})/2 +
\ell_1 \ell_2 (e^{2c} - 1)/2 \cr
&> a [{\rm Vol}(\partial V,k \vert_{\partial V}) +
\ell_1 \ell_2 (e^{2c} - 1)] /2 \cr
&= a {\rm Vol}(\partial V',k' \vert_{\partial V'})/2. \cr}$$
Therefore, for any $c>0$,
$\alpha(\ell_1 e^c) \geq \alpha(\ell_1)$, and hence
$\alpha$ is a non-decreasing function. $\square$

It is actually possible to define a function
$\alpha$ satisfying the conditions of Proposition 3.1
for which $\alpha(\ell_1) \rightarrow 1$
as $\ell_1 \rightarrow \infty$. In other words, if $\ell_1$ is
sufficiently large, then we can construct
a metric $k$ on $V$ for which the sectional
curvatures approach $-1$ and the volume approaches that of a cusp. 
This is intuitively plausible from the graphs of $f$ and $g$. 
However, a rigorous
proof of this result is slightly technical and long-winded.
Since we will not actually need this result,
we offer only a brief summary of the proof.

The idea is to construct, for any $t$
with $0 < t < 1$, the functions $f$ and $g$
in terms of a certain differential equation, which we omit
here. This differential equation in fact guarantees
that 
$$-1-t \leq \kappa_{\rm inf}(V,k)
\leq \kappa_{\rm sup}(V, k) \leq -1+t.$$
The constant $r_0$ is defined to be the
value of $r$ for which $f(r) = 0$. The
condition that $f'(r_0) = 2 \pi$ determines
$f(0)$, which is $\ell_1$. Hence, we obtain $\ell_1$
as a function of $t$. One shows that $\ell_1$
lies in the range $(2 \pi, \infty)$, and
that there exists an inverse function
$t \colon (2 \pi, \infty) \rightarrow (0,1).$
One also shows that  
as $\ell_1$ tends to $\infty$, the associated $t$ tends to zero.
In addition, the definition of the metric $k$ ensures that as
$\ell_1 \rightarrow \infty$, the ratio
$${{\rm Vol}(V,k) \over {\rm Vol}(\partial V, k \vert_{\partial V})/2}$$
tends to $1$. Hence, it is straightforward to construct
the function $\alpha$ satisfying the conditions of Proposition 3.1
and also $\alpha(\ell_1) \rightarrow 1$ as $\ell_1 \rightarrow \infty$.

We now apply Proposition 3.1 to $X(s_1, \dots, s_n)$, where $X$ is
a hyperbolic 3-manifold and $s_1, \dots, s_n$ are slopes on 
distinct components of
$\partial X$. The following proposition analyses the volume and
sectional curvatures of $X(s_1, \dots, s_n)$.

\noindent {\bf Proposition 3.2.} {\sl Let
$\alpha \colon (2\pi, \infty) \rightarrow (0,1)$ be the function
in Proposition 3.1. Let $X$ be a compact 3-manifold
with interior having a complete finite volume hyperbolic metric $h$,
and let $s_1, \dots, s_n$ be slopes on distinct tori
$T_1, \dots, T_n$ in $\partial X$.
Suppose that there is a maximal horoball neighbourhood of
$T_1 \cup \dots \cup T_n$ on which $l(s_i) > 2\pi$ for each $i$.
Let $\ell = \min_{1 \leq i \leq n} l(s_i)$.
Then $X(s_1, \dots, s_n)$ has a complete finite volume
negatively curved Riemannian metric $g$ for
which the following formulae hold.
$$-(\alpha(\ell))^{-1}
\leq \kappa_{\rm inf}(X(s_1, \dots, s_n),g),$$
$$\kappa_{\rm sup}(X(s_1, \dots, s_n),g) \leq
-\alpha(\ell),$$
$$\alpha(\ell) <
{{\rm Vol}(X(s_1, \dots, s_n),g) \over {\rm Vol}(X,h)}.$$}

\noindent {\sl Proof.} We follow the proof of the `$2 \pi$' theorem.
Let $\bigcup_{i=1}^n B_i$ be a union of horoballs in ${\Bbb H}^3$ which
projects to a
maximal horoball neighbourhood $N$ of $T_1 \cup \dots \cup T_n$.
Suppose that $B_i$ projects to the cusp at $T_i$.
Let $T'_i$ be the quotient of $\partial B_i$ by the subgroup
of parabolic isometries in $\pi_1(X)$ which preserve $B_i$. 
Now construct as in Proposition 3.1 a metric $k_i$ on the solid torus
$V_i$ which agrees on $\partial V_i$ with the Euclidean metric on $T'_i$,
and which has meridian length $l(s_i)$. The Riemannian
metric $g$ on $X(s_1, \dots, s_n)$ is just that obtained by
attaching $\bigcup_{i=1}^n (V_i,k_i)$ to $({\rm int}(X) - {\rm int}(N),h)$.

The first formula of Proposition 3.1 immediately
implies the first two formulae of the proposition.
To obtain the third formula, first note that it is elementary calculation
that
$${\rm Vol}(N,h) = \sum_{i=1}^n {\rm Vol}(T'_i, h \vert_{T'_i})/2.$$
Also, the metrics on $T'_i$ and $\partial V_i$ agree:
$${\rm Vol}(T'_i, h\vert_{T'_i})=
  {\rm Vol}(\partial V_i, k_i\vert_{\partial V_i}).$$
Proposition 3.1 gives that
$$\eqalign{{\rm Vol}(V_i, k_i) &\geq  \alpha(l(s_i))
{\rm Vol} (\partial V_i, k_i \vert_{\partial V_i}) /2 \cr
&\geq \alpha(\ell)
{\rm Vol} (\partial V_i, k_i \vert_{\partial V_i}) /2,\cr}$$
since $\alpha$ is a non-decreasing function. Hence,
$$\sum_{i=1}^n {\rm Vol}(V_i, k_i) \geq
\alpha(\ell)
\sum_{i=1}^n {\rm Vol} (\partial V_i, k_i \vert_{\partial V_i}) /2.$$
So,
$$\eqalign{
{\rm Vol}(X(s_1, \dots, s_n), g)
&= {\rm Vol}(X - N, h)+ \sum_{i=1}^n {\rm Vol}(V_i, k_i)\cr
&\geq {\rm Vol}(X - N, h) +
\alpha(\ell)
\sum_{i=1}^n {\rm Vol} (\partial V_i, k_i \vert_{\partial V_i}) /2 \cr
&> \alpha(\ell)
[{\rm Vol}(X - N, h) +
\sum_{i=1}^n {\rm Vol} (\partial V_i, k_i \vert_{\partial V_i}) /2 ]\cr
&= \alpha(\ell)
[{\rm Vol}(X - N, h) + {\rm Vol}(N, h)]\cr
&= \alpha(\ell)
[{\rm Vol}(X, h)],}$$
which establishes the final formula of the proposition. $\square$

The aim now is to use the comparison of volumes and sectional curvatures
of $(X, h)$ and $(X(s_1, \dots, s_n), g)$ to make a comparison of their
Gromov norms.

\noindent {\bf Proposition 3.3.} {\sl There is a non-increasing function
$\beta \colon (2 \pi, \infty) \rightarrow (1 , \infty)$, which has the
following property. Let $X$ be a compact hyperbolic 3-manifold
and let $s_1, \dots, s_n$ be slopes on distinct components
$T_1, \dots, T_n$ of $\partial X$. Suppose that there is
a maximal horoball neighbourhood of $T_1 \cup \dots \cup T_n$ on which
$l(s_i) > 2 \pi$ for each $i$. Then
$$\vert X(s_1, \dots, s_n) \vert \leq \vert X \vert
< \vert X(s_1, \dots, s_n) \vert \, \beta(\min_{1 \leq i \leq n} l(s_i)).$$
}

\noindent {\sl Proof.} It is a corollary of [6.5.2, 17]
that $\vert X(s_1, \dots, s_n) \vert \leq \vert X \vert$.
To prove the second half of the inequality of the proposition,
we need to compare the volume of a negatively curved Riemannian manifold
with its Gromov norm. For a manifold $X$ with interior having a complete
finite volume hyperbolic metric $h$, it is proved in [6.5.4, 17] that
$${\rm Vol}(X,h) = v_3 \, \vert X \vert,$$
where $v_3$ is the volume of a regular ideal 3-simplex in ${\Bbb H}^3$.
Here, $X(s_1, \dots, s_n)$ can, by the `$2 \pi$' theorem, be given
a negatively curved metric $g$, and so we need a version of
Gromov's result which applies in this case. It is possible to
show [C.5.8, 3] that if $M$ is a 3-manifold with a complete finite
volume Riemannian metric $g$
which has $\kappa_{\sup}(M,g) \leq -1$, then
$$ \vert M \vert \, \pi/2 \geq {\rm Vol}(M,g).$$
Now, if the metric $g$ on $X(s_1, \dots, s_n)$
is scaled by a positive constant $\lambda$
to give a metric $\lambda g$, then it is an elementary consequence
of the definition of sectional curvature that
$$\kappa_{\sup}(X(s_1, \dots, s_n), \lambda g) =
\lambda^{-2} \kappa_{\sup}(X(s_1, \dots, s_n), g).$$
Hence, by letting $\lambda = \sqrt{-\kappa_{\sup}(X(s_1, \dots, s_n), g)}$,
we deduce that
$$\eqalign{\vert X(s_1, \dots, s_n) \vert \pi /2
&\geq {\rm Vol}(X(s_1, \dots, s_n), \lambda g) \cr
&= \lambda^{3} {\rm Vol}(X(s_1, \dots, s_n), g)\cr
&=(-\kappa_{\sup}(X(s_1, \dots, s_n),g))^{3/2}
{\rm Vol}(X(s_1, \dots, s_n),g).}$$
We can now use Proposition 3.2 to deduce that
$$\eqalign{\kappa_{\sup}(X(s_1, \dots, s_n),g)
&\leq -\alpha(\ell)\cr
{\rm Vol}(X(s_1, \dots, s_n),g)
&> \alpha(\ell) [{\rm Vol}(X,h)],}$$
where $\alpha$ is the function given in Proposition 3.1
and $\ell = \min_{1 \leq i \leq n} l(s_i)$. So,
$$\eqalign{
\vert X(s_1, \dots, s_n) \vert \pi/2 &>
[\alpha(\ell)]^{5/2} {\rm Vol} (X,h)\cr
&= [\alpha(\ell)]^{5/2} v_3 \, \vert X \vert.}$$
The proposition is now proved by letting
$\beta(x) = [\alpha(x)]^{-5/2}\pi/2v_3$.
$\square$

It should be possible to find a function $\beta$
satisfying the requirements of Proposition 3.3 and for which
$\beta(\ell_1) \rightarrow 1$ as $\ell_1 \rightarrow \infty$.
To find such a $\beta$, one examines
the volume of a straight 3-simplex $\Delta$ in a simply-connected
3-manifold with sectional curvatures between
$-1 - \delta$ and $-1 + \delta$ for sufficiently small $\delta > 0$.
If one shows that as $\delta \rightarrow 0$, the
maximal volume of $\Delta$ tends to $v_3$, then one can find
a $\beta$ satisfying the conditions of Proposition 3.3, for which
$\beta(\ell_1) \rightarrow 1$ as $\ell_1 \rightarrow \infty$.

\vskip 18pt
\centerline {4. APPLICATIONS TO DEHN SURGERY}
\vskip 6pt

We now use the estimates of the previous section to deduce
some new results about Dehn surgery. The most far-reaching of these
is the following theorem.

\noindent {\bf Theorem 4.1.} {\sl Let $M$ be a compact
orientable 3-manifold, with
$\partial M$ a (possibly empty) union of tori.
Let $X$ be a hyperbolic manifold and let $s_1, \dots, s_n$
be a collection of slopes on $n$ distinct tori $T_1, \dots, T_n$
in $\partial X$,
such that $X(s_1, \dots, s_n)$ is homeomorphic to $M$.
Suppose that there exists in ${\rm int}(X)$ a maximal horoball neighbourhood
of $T_1 \cup \dots \cup T_n$
on which each slope $s_i$ has length at least
$2 \pi + \epsilon$, for some $\epsilon >0$.
Then, for any given $M$ and $\epsilon$, there is
only a finite number of possibilities (up to isometry)
for $X$, $n$ and $s_1, \dots, s_n$.}

To prove this, we shall need two well-known lemmas.

\noindent {\bf Lemma 4.2.} {\sl Let $X$ be a compact
orientable hyperbolic 3-manifold, and
let $T_1, \dots, T_n$ be a collection of tori in $\partial X$.
For each $i \in {\Bbb N}$ and $j \in \lbrace 1, \dots , n \rbrace$,
let $s_i^j$ be a slope on $T_j$.
Assume that, for every $j$, $s_i^j \not= s_k^j$ if $i \not= k$.
Then any given 3-manifold $M$
is homeomorphic to $X(s_i^1, \dots, s_i^n)$ for at most
finitely many $i$.}

\noindent {\sl Proof.} By the hyperbolic Dehn surgery theorem of Thurston
[Theorem E.5.1, 3],
$X(s_i^1, \dots, s_i^n)$ is hyperbolic for $i$ sufficiently large.
Thus, if the theorem were not true, then $M$ would have to be hyperbolic.
Moreover, from the proof of Thurston's hyperbolic Dehn
surgery theorem, for $i$ sufficiently large, the cores of the
filled-in solid tori in $X(s_i^1, \dots, s_i^n)$ are geodesics,
whose lengths each tend to zero, as $i \rightarrow \infty$.
In particular, for $i$ sufficiently large,
$X(s_i^1, \dots, s_i^n)$ has a geodesic shorter
than the shortest geodesic in $M$. This is impossible by
Mostow rigidity [Theorem C.5.4, 3]. $\square$

\noindent {\bf Lemma 4.3.} {\sl Let $X$ be a compact
orientable hyperbolic 3-manifold, and
let $T_1, \dots, T_n$ be a collection of tori in $\partial X$.
Let $\lbrace (s_i^1, \dots, s_i^n) : i \in {\Bbb N} \rbrace$
be a sequence of distinct $n$-tuples, where each $s_i^j$ is a slope on $T_j$.
Suppose that, for each $i$, we can find a maximal horoball neighbourhood
of the cusps at $T_1 \cup \dots \cup T_n$ on which each $s_i^j$ has length
more than $2 \pi$.
Then any given 3-manifold $M$
is homeomorphic to $X(s_i^1, \dots, s_i^n)$ for at most
finitely many $i$.}

\noindent {\sl Proof.}
If the lemma were not true, we could pass to a subsequence, such that
$M$ is homeomorphic to $X(s_i^1, \dots, s_i^n)$ for each $i$.
There are two possibilities for each $j \in
\lbrace 1, \dots, n \rbrace$:

\item{(i)} the sequence $\lbrace s_i^j : i \in {\Bbb N} \rbrace$
contains a subsequence in which the slopes $s_i^j$ are all distinct, or
\item{(ii)} the sequence $\lbrace s_i^j : i \in {\Bbb N} \rbrace$
runs through only finitely many slopes.

Since the $n$-tuples $(s^1_i, \dots, s^n_i)$ are distinct,
at least one $j$ satisfies (i). After re-ordering,
we may assume this value of $j$ is $n$.
Pass to this subsequence, where the slopes $s^n_i$ are all distinct.
In this new sequence, the integers $j \in \lbrace 1, \dots, n-1 \rbrace$
either satisfy (i) or (ii). If some $j$ satisfies (i), say $j = n-1$,
pass to this subsequence.
Continuing in this fashion, we obtain a sequence and
an integer $m \geq 1$ such that
\item{(i)} $s_i^j \not= s^j_k$ for $i \not= k$ and $j \geq m$, and
\item{(ii)} $\lbrace s_i^j : i \in {\Bbb N} \rbrace$ runs through
only finitely many slopes, for each $j < m$.

By passing to a subsequence, we may assume that $s_i^j$ is
the same slope $s^j$ for all $i$, when $j < m$.
If $m > 1$, let $Y = X(s^1, \dots, s^{m-1})$. Otherwise,
let $Y = X$. Now, we may find a maximal horoball neighbourhood of the
cusps at $T_1 \cup \dots \cup T_n$ on which each $s^j$ has length
more than $2 \pi$. Hence, $Y$ admits
a negatively curved metric, by the `$2 \pi$' theorem.
Hence, it cannot be reducible, toroidal or Seifert fibred.
(See the proof of Proposition 2.2.)
Its boundary contains $T_m \cup \dots \cup T_n$ and so is non-empty.
Hence, by Thurston's theorem
on the geometrisation of Haken 3-manifolds [Chapter V, 12], $Y$ is hyperbolic.
But, $M$ is homeomorphic to $Y(s_i^m, \dots, s_i^n)$ for each $i$.
Lemma 4.2 gives us a contradiction. $\square$

\noindent {\sl Proof of Theorem 4.1.}
Suppose that, on the contrary, there exists
a sequence of 3-manifolds $X_i$ with complete finite volume
hyperbolic metrics $h_i$ on their interiors,
and slopes $s^1_i, \dots, s^{n(i)}_i$ on $\partial X_i$
with $l(s^j_i) \geq 2 \pi + \epsilon$,
such that $X_i(s^1_i, \dots, s^{n(i)}_i)$
is homeomorphic to $M$. Then, by Proposition 3.3,
$$\vert X_i \vert < \vert M \vert \,
\beta(\min_{1 \leq j \leq n(i)} l(s_i^j)).$$
Since $\beta$ is a non-increasing function,
$$\vert X_i \vert < \vert M \vert \, \beta(2 \pi + \epsilon).$$
Thus, the sequence $\vert X_i \vert$ is bounded, and so
the sequence ${\rm Vol}(X_i,h_i)$ is also bounded,
since the Gromov norm and the volume of a hyperbolic 3-manifold are
proportional [6.5.4, 17]. But for any real number $c$,
the collection of complete orientable hyperbolic 3-manifolds with volume
at most $c$ is a compact topological space when
endowed with the geometric topology [Theorem E.1.10, 3].
Hence, we may pass to a subsequence (also denoted $\lbrace X_i \rbrace$),
such that ${\rm int}(X_i)$ converges in the geometric topology to a
complete finite volume hyperbolic 3-manifold ${\rm int}(X_\infty)$, say,
where $X_\infty$ is compact and orientable.
This implies (see [Theorem E.2.4, 3]) that, for $i$
sufficiently large, the following is true.
In each 3-manifold ${\rm int}(X_i)$, there is a (possibly empty) union $L_i$
of disjoint closed geodesics, such that ${\rm int}(X_i) - L_i$
is diffeomorphic to ${\rm int}(X_\infty)$. This
diffeomorphism is a $k_i$-bi-Lipschitz map except in
a small neighbourhood of $L_i$, for real numbers $k_i \geq 1$
which tend to $1$, as $i \rightarrow \infty$.
This diffeomorphism also takes a maximal horoball neighbourhood $N_i$
of cusps of ${\rm int}(X_i)$ to a neighbourhood of cusps
of ${\rm int}(X_\infty)$ which closely approximates a
maximal horoball neighbourhood $N'_i$.
We extend $N'_i$ to a maximal horoball neighbourhood $N_i''$
of all the cusps of ${\rm int}(X_\infty)$.
The slopes $s_i^1, \dots, s_i^{n(i)}$ correspond to slopes
$\sigma_i^1, \dots, \sigma_i^{n(i)}$
on toral components of $\partial X_\infty$, and
the ratios $l(s_i^j)/l(\sigma_i^j) \rightarrow 1$
as $i \rightarrow \infty$, where the lengths $l(s_i^j)$
and $l(\sigma_i^j)$ are measured on $N_i$ and $N''_i$
respectively. Hence, as $l(s_i^j) \geq 2 \pi + \epsilon$,
$l(\sigma_i^j) > 2 \pi$ for $i$ sufficiently large.
Since $X_\infty$ has a finite number of cusps,
the sequences $n(i)$ and $\vert L_i \vert$ are bounded. Hence, by passing to
a subsequence, we may assume that $n(i)$ is some fixed positive
integer $n$ and that $\vert L_i \vert$ is some fixed non-negative
integer $p$, for all $i$.
We may pass to a further subsequence
where for any $j$, $\sigma_i^j$ lies on the same torus
for all $i$. Now, $p > 0$, for otherwise, $X_\infty = X_i$ for all $i$
and then $M$ is homeomorphic to $X_\infty(\sigma_i^1, \dots, \sigma_i^{n})$
for each $i$. This is a contradiction by Lemma 4.3.
Thus there are slopes $(t^1_i, \dots, t^p_i)$ on
$\partial X_\infty$ such that $X_\infty(t^1_i, \dots, t^p_i)$
is homeomorphic to $X_i$. We may assume that, for any $j$,
the slopes $t^j_i$ lie on the same torus $P^j$ for all $i$.
Since the manifolds ${\rm int}(X_i)$ converge in the geometric topology to
${\rm int}(X_\infty)$, we may assume that the slopes $t^j_i$
are all distinct. We shall now show that, for each $j$,
$l(t^j_i) \rightarrow \infty$, as $i \rightarrow \infty$,
where the slope lengths are measured on $N_i''$.
For each torus $T^k$ in $\partial X_\infty$, let $N^k$ be
the maximal horoball neighbourhood of $T^k$.
Then, we may find a horoball neighbourhood $H^j$ of $P^j$ which
misses all $N^k$ other than the horoball $N^j$ corresponding
to $P^j$. Thus, $H^j$ lies inside $N''_i$ for each $i$.
Now, the lengths of $t^j_i$ tend to $\infty$ as $i \rightarrow
\infty$, where the length is measured on $H^j$, since
the slopes $t^j_i$ are all distinct.
Hence, the lengths $l(t^j_i) \rightarrow \infty$, as $i \rightarrow \infty$,
where the slope lengths are measured on $N_i''$.
So, $M$ is homeomorphic to
$X_\infty(t^1_i, \dots, t^p_i, \sigma_i^1, \dots, \sigma_i^n)$ for
each $i$, and for sufficiently large $i$,
$l(t^j_i) > 2 \pi$ and $l(\sigma_i^j) > 2 \pi$
for each $j$, where the slope lengths are measured
on $N''_i$. This is a contradiction, by Lemma 4.3. $\square$

Theorem 4.1 has the following corollary.

\noindent {\bf Theorem 4.4.} {\sl Let $M$ be a compact orientable 3-manifold
with $\partial M$ a (possibly empty) union of tori. Suppose that $M$ is
homeomorphic to $X(s)$, where $X$ is a hyperbolic 3-manifold
and $s$ is a slope on a toral boundary component $T$ of $X$.
Suppose also that $e$ is a short slope on $T$
such that $\Delta(s,e) > 22$, or that
$e$ is a minimal slope with $\Delta(s,e) > 3$.
Then, for a given $M$, there is
only a finite number of possibilities (up to isometry) for $X$, $s$ and $e$.}

\noindent {\sl Proof.}
Suppose that $e$ is a short slope with $\Delta(s,e) > 22$.
The proof when $e$ is minimal is entirely analagous.
Fix $\epsilon$ as $(23 \sqrt 3/2 \pi) - 2\pi$,
which is greater than zero.
By Theorem 4.1, there is only a finite number of hyperbolic 3-manifolds
$X$ and slopes $s$ on a torus component of $\partial X$, such that
$X(s)$ is homeomorphic to $M$ and such that $l(s) \geq 2 \pi + \epsilon$. But
if $e$ is a short slope on $T$ with $\Delta(s,e) \geq 23$,
then, by Corollary 2.4,
$$l(s) \geq {\sqrt 3 \, \Delta (s,e) \over 2 \pi} \geq
{23 \sqrt 3 \over 2 \pi} = 2\pi + \epsilon.$$
Thus, there is only a finite number of possibilities
for $X$ and $s$. Also there is only a finite number
of short slopes $e$ on $T$. Hence, the theorem is proved. $\square$

\noindent {\bf Corollary 4.5.} {\sl For a given closed orientable
3-manifold $M$, there is at most a finite number of hyperbolic
knots $K$ in $S^3$
and fractions $p/q$ (in their lowest terms) such that $M$ is
obtained by $p/q$-Dehn surgery along $K$ and $\vert q \vert > 22$.}

\noindent {\sl Proof.} The meridian slope $e$ on
$\partial {\cal N}(K)$ is a short slope, and
$\Delta(e, p/q) = \vert q \vert$. Now apply Theorem 4.4. $\square$

\noindent {\bf Theorem 4.6.} {\sl Let $M_1$ and $M_2$ be
compact orientable 3-manifolds with $\partial M_i$ a (possibly empty)
union of tori. Let $X$ be a hyperbolic 3-manifold and
let $T$ be toral boundary component of $X$. Suppose that there
are slopes $s_1$ and $s_2$ on $T$, with $\Delta (s_1, s_2) > 22$,
such that $X(s_i)$ is homeomorphic to $M_i$ for $i =1$ and $2$.
Then, for any given $M_1$ and $M_2$,
there is only a finite number of possibilities (up to isometry) for
$X$, $s_1$ and $s_2$.}

\noindent {\sl Proof.}
Suppose that there is an infinite number of triples $(X, s_1, s_2)$,
for which $X(s_i)$ is homeomorphic to $M_i$ (for $i =1$ and $2$)
and which have $\Delta(s_1, s_2) \geq 23$.
Let $\epsilon = (\root 4 \of 3 \, \sqrt {23}) - 2\pi$.
If $l(s_1) < 2 \pi + \epsilon$ and $l(s_2) < 2 \pi + \epsilon$, then
$$l(s_1) \, l(s_2) <  (2 \pi + \epsilon)^2 = 23 \sqrt 3 \leq
\Delta(s_1,s_2) \sqrt 3,$$
but this cannot occur, by Lemma 2.1. Hence, an
infinite number of the slopes $s_1$ or $s_2$ must have
have length at least $2 \pi + \epsilon$.
For the sake of definiteness, assume $l(s_1) \geq 2\pi + \epsilon$,
for an infinite number of slopes $s_1$.
By Theorem 4.1, there is only a finite number of possibilities for $X$
and $s_1$. Hence, by passing to an infinite subcollection,
we can find a fixed hyperbolic manifold $X$ and an infinite number
of distinct slopes $s_2$ such that $X(s_2)$ is homeomorphic to $M_2$.
This contradicts Lemma 4.2. $\square$

The next result is a `uniqueness' theorem for Dehn surgery.
It should be compared with [Theorem 4.1, 11].

\noindent {\bf Theorem 4.7.} {\sl For $i = 1$ and $2$, let
$X_i$ be a hyperbolic 3-manifold and let $s_i$ be a slope on
$\partial X_i$. Then, there is a
real number $C(X_1)$ depending only on $X_1$,
such that if $l(s_2) > C(X_1)$, then
$$\left\{ X_1(s_1) \cong X_2(s_2) \right\}
\Longleftrightarrow \left\{ (X_1, s_1) \cong (X_2, s_2) \right\},$$
where $\cong$ denotes a homeomorphism.}

\noindent {\sl Proof.} Suppose that there is no such real number
$C(X_1)$. Then, there is a sequence of hyperbolic 3-manifolds $X^i_2$
and slopes $s^i_2$ on $\partial X^i_2$
such that $l(s^i_2) \rightarrow \infty$,
together with slopes $s^i_1$ on $\partial X_1$,
with the property that $X_1(s^i_1) \cong X^i_2(s^i_2)$, but
$(X_1, s^i_1) \not\cong (X^i_2, s^i_2)$.

Note first that the sequence $s^i_1$ can have no constant subsequence.
For, if there were such a subsequence, say with slope $s_1$,
then $X_1(s_1) \cong X^i_2(s^i_2)$ for infinitely many $i$.
This contradicts Theorem 4.1.

\noindent {\sl Case 1.} $X^i_2$ runs through only finitely many
hyperbolic 3-manifolds (up to isometry).

Then, by passing to a subsequence, we may assume that $X^i_2$ is
some fixed hyperbolic manifold $X_2$.
Since $s^i_1$ cannot run through only finitely many slopes,
we may pass to a subsequence where
$s^i_1 \not= s^j_1$ if $i \not= j$.
By Thurston's hyperbolic Dehn surgery theorem,
for $i$ sufficiently large, $X_1(s^i_1)$ is hyperbolic, with
the core of the surgery solid torus being the unique shortest
geodesic. Similarly, for $i$ sufficiently large,
$X_2(s^i_2)$ is hyperbolic, with
the core of the surgery solid torus being the unique shortest
geodesic. But, then by Mostow Rigidity, there is a homeomorphism
from $X_1(s^i_1)$ to $X_2(s^i_2)$ which takes one geodesic to the
other. Hence, $(X_1, s^i_1) \cong (X_2, s^i_2)$, which is contrary
to assumption. Hence, the following case must hold.

\noindent {\sl Case 2.} There is a subsequence in which
$X^i_2$ and $X^j_2$ are not isometric if $i \not= j$.

Now, by Proposition 3.3,
$$(\beta (l(s^i_2)))^{-1} \, \vert X^i_2 \vert
< \vert X^i_2(s^i_2) \vert = \vert X_1(s^i_1) \vert
\leq \vert X_1 \vert.$$
Since $\beta$ is a non-increasing function,
the sequence $\vert X^i_2 \vert$ is bounded, and so,
by passing to a subsequence, we may assume that the
3-manifolds ${\rm int}(X^i_2)$ converge in the geometric topology
to a hyperbolic 3-manifold ${\rm int}(X^\infty_2)$,
where $X^\infty_2$ is compact and orientable. 
For $i$ sufficiently large, there is a union $L_i$ of $n>0$
disjoint closed geodesics in ${\rm int}(X^i_2)$,
such that ${\rm int}(X^i_2) - L_i$
is diffeomorphic to ${\rm int}(X^\infty_2)$.
In the complement of a small neighbourhood of $L_i$,
this map is $k_i$-Lipschitz
for real numbers $k_i \geq 1$ which tend
to $1$. Let $(t^i_1, \dots, t^i_n)$ be the slopes
on $\partial X^\infty_2$ such that $X^\infty_2(t^i_1, \dots, t^i_n)$
is homeomorphic to $X^i_2$. By passing to a subsequence, we may
ensure that the slopes $t^i_j$ are all distinct.
Now, the slopes $s^i_2$ correspond to slopes $\sigma^i_2$, say
on $\partial X^\infty_2$. Since the lengths $l(s^i_2) \rightarrow \infty$,
so also the lengths $l(\sigma^i_2) \rightarrow \infty$.
Hence, by Thurston's hyperbolic Dehn surgery theorem,
for any $\epsilon > 0$, $X^\infty_2(\sigma^i_2, t^i_1, \dots, t^i_n)$
is hyperbolic and the cores of the $(n+1)$ filled-in solid tori
are geodesics of length less than $\epsilon$, if $i$ is
sufficiently large.
However, $X^\infty_2(\sigma^i_2, t^i_1, \dots, t^i_n)$ is
homeomorphic to $X_1(s^i_1)$.
Since $s^i_1$ has no constant subsequence, Thurston's
hyperbolic Dehn surgery theorem gives that there is an
integer $N$ and an $\epsilon >0$ such that, for all $i \geq N$,
$X_1(s^i_1)$ is hyperbolic and
the core of the filled-in solid torus is the
unique geodesic with length less than $\epsilon$.
This is a contradiction. $\square$

Theorem 4.1 also has the following corollary regarding
branched covers.

\noindent {\bf Corollary 4.8.} {\sl Let $M$ be a compact
orientable 3-manifold with $\partial M$ a (possibly empty)
union of tori, which is obtained as a branched
cover of a compact orientable 3-manifold $Y$ over a hyperbolic link $L$,
via representation $\rho \colon \pi_1(Y-L) \rightarrow S_r$.
Suppose that the branching index of every lift
of every component of $\partial {\cal N}(L)$ is
at least 7. Then, for a given $M$, there are only finitely
many possibilities for $Y$, $L$, $r$ and $\rho$.}

\noindent {\sl Proof.} The representation $\rho \colon \pi_1(Y - L)
\rightarrow S_r$ determines a cover $X$ of $Y - {\rm int}({\cal N}(L))$,
and $M$ is obtained from $X$ by Dehn filling along slopes
$s_1, \dots, s_n$ in $\partial X$. The hyperbolic structure
on $Y - {\rm int}({\cal N}(L))$ lifts to a hyperbolic
structure on $X$, and a maximal horoball neighbourhood
of $\partial {\cal N}(L)$ lifts to a horoball neighbourhood
of cusps of $X$. Since the length of each slope
on $\partial {\cal N}(L)$ in $Y - {\rm int}({\cal N}(L))$
is at least 1 [4], the length of each $s_i$ on $N$ is at
least 7. Thus Theorem 4.1 implies that there are only
finitely many possibilities for $X$ and $s_1, \dots s_n$.
Now,
$${\rm Vol}(X) = r {\rm Vol}(Y - {\rm int}({\cal N}(L))),$$
where $r$ is the index of the cover $X \rightarrow 
Y - {\rm int}({\cal N}(L))$. There is a lower bound on ${\rm Vol}
(Y - {\rm int}({\cal N}(L)))$, since the volume of
a maximal horoball neighbourhood of $\partial {\cal N}(L)$
is at least $\sqrt 3$ [1]. Hence, for a given $X$, there
is an upper bound on $r$. Once $r$ and $X$ are fixed,
so is ${\rm Vol}(Y - {\rm int}({\cal N}(L)))$. 
There are only finitely many hyperbolic manifolds of
a given volume [3], and so there are only finitely many
possibilities for $Y - {\rm int}({\cal N}(L))$. The
lengths of the meridian slopes on $\partial {\cal N}(L)$
are bounded above by the lengths of the slopes
$s_i$ on $N$.
Hence, there are only finitely many possibilities for
$Y$ and $L$. A representation $\rho \colon \pi_1(Y - L)
\rightarrow S_r$ is determined by the image of a
generating set of $\pi_1(Y - L)$. Hence, once $Y$, $L$
and $r$ are fixed, there are only finitely many possibilities for $\rho$.
$\square$

\vfill\eject
\centerline {5. THE LENGTH OF BOUNDARY SLOPES}
\vskip 6pt

The results about Dehn surgery in Section 4 bear a strong resemblance
to the work in [11]. In that paper, the main theorem (1.4) of [10]
was crucial in establishing strong restrictions on
the number of intersection points between embedded surfaces in
a 3-manifold and certain surgery curves. In this section, we
use hyperbolic techniques to prove similar results.
The main theorem of this section is the following.

\noindent {\bf Theorem 5.1.} {\sl Let $X$ be a hyperbolic 3-manifold
and let $T$ be a toral component of $\partial X$. Let
$f \colon F \rightarrow X$ be a map of a compact connected surface $F$
into $X$, such that $f(F) \cap \partial X = f(\partial F)$.
Suppose that $f_\ast \colon \pi_1(F) \rightarrow \pi_1(X)$
is injective, and that every essential arc in $F$ maps to an
arc in $X$ which cannot be homotoped (rel its endpoints) into $\partial X$.
Suppose also that $f(F) \cap T$ is a non-empty
collection of disjoint simple closed curves, each with
slope $s$. Then
$$l(s) \, \vert f(F) \cap T \vert < -2 \pi \, \chi(F).$$}

This result has a number of corollaries, which include the following.

\noindent {\bf Corollary 5.2.} {\sl Let $p/q$ be the boundary slope of 
an incompressible boundary-incompressible non-planar orientable
surface $F$ properly embedded in the exterior of a hyperbolic knot
in $S^3$. Then
$$\vert q \vert < {4 \pi^2  \, {\rm genus}(F) \over \sqrt 3}.$$}

\noindent {\sl Proof.} By Theorem 5.1,
$$l(p/q) \, \vert \partial F \vert < - 2 \pi \, \chi(F)
=2 \pi (2 \, {\rm genus}(F) -2 + \vert \partial F \vert),$$
and so
$$(l(p/q) - 2 \pi) \vert \partial F \vert < 4 \pi ({\rm genus}(F) - 1).$$
Now, the inequality of the corollary
is automatically satisfied when $\vert q \vert = 1$,
since ${\rm genus}(F) > 0$.
Hence we may assume that $\vert \partial F \vert \geq 2$.
Hence, if $l(p/q) \geq 2\pi$,
$$2 l(p/q) - 4 \pi \leq (l(p/q) - 2 \pi) \vert \partial F \vert
< 4 \pi ({\rm genus}(F) - 1).$$
If $l(p/q) < 2 \pi$, then
$$2 l(p/q) - 4 \pi < 0 \leq 4 \pi ({\rm genus}(F) - 1).$$
So, in either case,
$$l(p/q) < 2 \pi \, {\rm genus}(F).$$
By Proposition 2.2, the meridian slope on
$\partial {\cal N}(K)$ is short. By Corollary 2.4,
$$l(p/q) \geq \vert q \vert \, \sqrt 3 / 2 \pi.$$
Hence, we obtain the inequality of the corollary. $\square$

Let $F$ be as in Theorem 5.1. Then $F$ is neither a disc, nor a
M\"obius band, nor an annulus. Since
$[s] \in \pi_1(X)$ is non-trivial,
$F$ cannot be a disc. If $F$ were a
M\"obius band or an annulus, then the map $f$
could be homotoped (rel $\partial F$) to a
map into $\partial X$, as $X$ is hyperbolic.
Since $F$ contains an essential arc,
this is contrary to assumption.

We may therefore pick an ideal triangulation of ${\rm int}(F)$. 
In other words, we may express
${\rm int}(F)$ as a union of 2-simplices glued along
their edges, with the 0-simplices then removed.
We may also ensure that each 1-cell in $F$ is an essential arc.
To see that such an ideal triangulation exists,
fill in the boundary components of $F$ with discs,
forming a closed surface $F^+$. If $F^+$ has non-positive
Euler characteristic, then it admits a one-vertex triangulation,
in which each 1-cell is essential. (By a `triangulation'
here, all we mean is an expression of $F^+$ as union of 2-simplices
with their edges identified in pairs.) If $F^+$
is a projective plane, then it admits a two-vertex
triangulation. If $F^+$ is a sphere,
then it admits a triangulation with three vertices.
After subdividing, if necessary, each of these triangulations
to increase the number of vertices, 
we obtain an ideal triangulation of $F$ of the required
form.

The following result is due to Thurston [Section 8, 17].

\noindent {\bf Proposition 5.3.} {\sl There is a homotopy of 
$f \colon F \rightarrow X$
to a map which sends each ideal triangle of ${\rm int}(F)$ 
to a totally geodesic ideal triangle in ${\rm int}(X)$.}

\noindent {\sl Proof.}
We construct the homotopy on the 1-cells first.
First pick a horoball neighbourhood $N$ of the cusps of ${\rm int}(X)$
which is a union of disjoint copies of $S^1 \times S^1 \times [1, \infty)$.
$N$ lifts to a disjoint union of horoballs in ${\Bbb H}^3$.
We may homotope $f$ so that, after the
homotopy, $f(F) \cap N$ is a union of vertical half-open
annuli. Hence, for each (open) 1-cell $\alpha$ in
the ideal triangulation, $\alpha - f^{-1}({\rm int}(N))$ is a single interval.
We may homotope this interval, keeping its endpoints fixed, to a
geodesic in ${\rm int}(X)$.  By assumption, this geodesic does not
wholly lie in $N$. Hence, in the universal cover,
this geodesic runs between distinct horoball lifts of $N$. We can
then perform a further homotopy so that the whole of $\alpha$ is sent to
a geodesic.

The boundary of each 2-cell of $F^+$ is a union of three 1-cells.
The interior of each 1-cell is sent to a geodesic in ${\rm int}(X)$. By examining
the universal cover of ${\rm int}(X)$, it
is clear that we may map this 2-cell to an ideal triangle.
Furthermore, the map of this 2-cell is homotopic to the 
original map, since the universal cover of ${\rm int}(X)$ is aspherical.
$\square$

When $F$ is in this form, it is an example of a `pleated surface'.
It inherits a metric,
by pulling back the metric on ${\rm int}(X)$. This in fact gives
${\rm int}(F)$ a hyperbolic structure, since
the metric arises from glueing a union
of hyperbolic ideal triangles along their geodesic boundaries.
Furthermore, this structure is complete, since the metric on
${\rm int}(X)$ is complete.

\noindent {\sl Proof of Theorem 5.1.}
We may use Proposition 5.3 to homotope $f$ to a map $g$
such that $g({\rm int}(F))$ is a union of ideal triangles.
Let $N$ be the maximal horoball neighbourhood of $T$.
Let $N_-$ be a slightly smaller
horoball neighbourhood of $T$, such that $g^{-1}
(\partial N_-)$ is a disjoint union of simple closed curves,
and so that $g(F)$ intersects $\partial N_-$ transversely.
We may find a sequence of such $N_-$ whose
union is the interior of $N$.
We may also find a horoball neighbourhood $N'$ of $T$, strictly
contained in $N_-$,
such that $g(F) \cap N'$ is a union of vertical half-open annuli.
We will examine the intersection of $g(F)$ with
the region $N_- -{\rm int}(N')$, which is a copy
of $T^2 \times I$.
Consider a component $Y$ of $g^{-1} (N_- - {\rm int}(N'))$
for which $g(Y) \cap \partial N'$ is non-empty.

\noindent {\sl Claim 1.} $\partial Y$
contains no curve which bounds a disc in $F$.

Let $C$ be such a curve, bounding a disc $D$ in $F$.
Then $\partial D$ cannot map to $\partial N'$,
since $[s] \in \pi_1(X)$ is non-trivial. Thus,
the interior of $D$ is disjoint from $Y$.
If $C$ intersected no 1-cells of $F$, then $D$ would map
into $X$ in a totally geodesic fashion. Hence,
$g(D)$ would lie in $N_-$, and so $D$ would be $Y$.
This is impossible, and so
$C$ must intersect some 1-cells of $F$.
Hence, there is an arc in a 1-cell of $F$ which is
embedded $D$. But this 1-cell of $F$ maps to
a geodesic in $X$. This geodesic lifts to a geodesic in ${\Bbb H}^3$
which leaves and re-enters the same horoball lift of $N_-$. This
cannot happen.

\noindent {\sl Claim 2.} $g_\ast \colon \pi_1(Y) \rightarrow \pi_1(N_- -
{\rm int}(N'))$ is injective.

Suppose that $x$ is a non-zero element of $\pi_1(Y)$
which maps to $0 \in \pi_1(N - {\rm int}(N'))$.
Then $g_\ast \colon \pi_1(Y) \rightarrow \pi_1(X)$
sends $x$ to $0$. But $\pi_1(Y) \rightarrow \pi_1(X)$
factors as $\pi_1(Y) \rightarrow \pi_1(F) \rightarrow \pi_1(X)$.
The second of these maps is assumed to be injective.
Therefore $\pi_1(Y) \rightarrow \pi_1(F)$ sends $x$ to $0$.
This contradicts Claim 1.

\noindent {\sl Claim 3.}
$Y$ is an annulus with one boundary component mapping to
$\partial N'$ and the other mapping to $\partial N_-$.

Now, $\pi_1(N_- - {\rm int}(N')) \cong {\Bbb Z} \oplus {\Bbb Z}$.
Hence, by Claim 2, $Y$ must be a disc, a M\"obius band or an annulus.
However, if $Y$ is not an annulus with one boundary component
mapping to $\partial N'$ and the other mapping to $\partial N_-$,
then $F$ is a disc, a M\"obius band or an annulus,
which is a contradiction.

\noindent {\sl Claim 4.} Each component of $g^{-1}({\rm int}(N))$
which touches $f^{-1}(T)$ is an open annulus in ${\rm int}(F)$.

There is a sequence of horoball neighbourhoods $N_-$ whose
union is the interior of $N$. For each such $N_-$ we showed
in Claim 3 that each component $Y$ of $g^{-1}(N_-)$
which touches $f^{-1}(T)$ is a half-open annulus
in ${\rm int}(F)$. The union of these half-open annuli
is the required collection of open annuli in ${\rm int}(F)$.

There is a standard homeomorphism which identifies
$N$ with $S^1 \times S^1 \times [1, \infty)$.
We may pick such an identification so that
$S^1 \times \lbrace \ast \rbrace \times \lbrace 1 \rbrace$
has slope $s$, where $\lbrace \ast \rbrace$ is some point
in $S^1$. Consider now the covering
$$S^1 \times {\Bbb R} \times (1, \infty) \rightarrow
S^1 \times S^1 \times (1, \infty) \cong
{\rm int}(N)$$
which is determined by the subgroup generated by $[s] \in \pi_1(T)$.
Each open annulus from Claim 4 lifts to an open annulus
in $S^1 \times {\Bbb R} \times (1, \infty)$.
Now, $S^1 \times \lbrace 0 \rbrace \times (1, \infty)$
inherits a hyperbolic structure from ${\rm int}(N)$,
which makes it isometric to a 2-dimensional horocusp.
Let the map
$p\colon S^1 \times {\Bbb R} \times (1, \infty) \rightarrow
S^1 \times \lbrace 0 \rbrace \times (1, \infty)$
be orthogonal projection onto this submanifold.
Note that $p$ need not respect the product structure
of $S^1 \times {\Bbb R} \times (1, \infty)$.
Each open annulus of Claim 4 is mapped surjectively onto
$S^1 \times \lbrace 0 \rbrace \times (1, \infty)$.
Also, since $p$ does not increase distances,
the area $A'$ that each open annulus inherits from
$S^1 \times \lbrace 0 \rbrace \times (1, \infty)$
is no more than the area which it inherits from ${\rm int}(N)$.
However, $A'$ is at least $l(s)$, since this is the
area of the 2-dimensional horocusp
$S^1 \times \lbrace 0 \rbrace \times (1, \infty)$.
Thus, the hyperbolic area of $F$ is more than
$l(s) \, \vert f(F) \cap T \vert$.
But the Gauss-Bonnet formula [Proposition B.3.3, 3] states that
its total area is $-2 \pi \chi(F)$. We therefore deduce that
$$l(s) \, \vert f(F) \cap T \vert < -2 \pi \, \chi(F). \eqno \square$$

\vskip 18pt
\centerline {6. `ALMOST HYPERBOLIC' 3-MANIFOLDS ARE HYPERBOLIC}
\vskip 6pt

Throughout this paper, we have studied 3-manifolds which have
a complete finite volume negatively curved Riemannian metric.
It is a major conjecture whether the existence
of such a metric on a 3-manifold implies the existence of a
complete finite volume hyperbolic structure. In this
section, we provide evidence for this conjecture by
considering 3-manifolds which are `almost
hyperbolic' in the following sense.

\noindent {\bf Definition 6.1.} Let $\delta$
be a positive real number. Let $M$ be a compact orientable 3-manifold
with $\partial M$ a (possibly empty) collection of tori.
Let $g$ be a Riemannian metric on ${\rm int}(M)$.
Then $(M,g)$ is {\sl $\delta$-pinched} if
$$-1 - \delta \leq
\kappa_{\inf}(M,g) \leq \kappa_{\sup}(M,g) \leq -1 + \delta.$$
We say that $M$ is {\sl almost hyperbolic} if, for all
$\delta>0$, there a $\delta$-pinched complete
finite volume Riemannian metric on its interior.

The main theorem of this section is the following result.

\noindent {\bf Theorem 6.2.} {\sl Let $M$ be a compact orientable
3-manifold with $\partial M$ a (possibly empty) union of tori.
If $M$ is almost hyperbolic, then it has a complete finite volume
hyperbolic structure.}

We proved this theorem in the course of proving several
other results in this paper. We recently discovered that it has also
been proved by Zhou [18] using methods similar to our own (but
not identical).
It was also known to Petersen [14]. However, since the journal [18] is
relatively poorly circulated in the West, it seems worthwhile
to include a summary of our proof of this result here.

The idea behind our proof of this theorem is as follows.
Since $M$ is almost hyperbolic, there is
a sequence of positive real numbers $\delta_i$
tending to zero, and complete finite volume Riemannian
metrics $g_i$ on ${\rm int}(M)$ such that $(M, g_i)$ is $\delta_i$-pinched.
We show that some subsequence `converges' to a `limit'
manifold $(M_\infty, g_\infty)$ which is a 3-manifold $M_\infty$
with a complete finite volume hyperbolic hyperbolic metric $g_\infty$.
Hence, $M_\infty$ is the interior of some compact orientable
3-manifold $\overline M_\infty$ with $\partial \overline M_\infty$
a (possibly empty) union of tori. If $\overline M_\infty$
is homeomorphic to $M$, then we have found a complete
finite volume hyperbolic structure on $M$. There is no immediate
reason why $\overline M_\infty$ should be homeomorphic
to $M$, but we will show that, if it is not, then there exist slopes
$(s^1_i, \dots, s^{n(i)}_i)$ on $\partial \overline M_\infty$ such that
$\overline M_\infty(s^1_i, \dots, s^{n(i)}_i)$ is
homeomorphic to $M$. The length of each of these slopes tends
to infinity as $i \rightarrow \infty$. Lemma 4.2
then gives us a contradiction.

In Section 4, we exploited the well-known theory of
convergent sequences of hyperbolic manifolds, and in that
case, non-trivial convergence corresponds to hyperbolic
Dehn surgery [3]. In the case here, the infinite sequence of
manifolds are not hyperbolic, merely negatively curved, but
a similar theory applies. We recall the following definition [8],
due to Gromov (see also [6]).

\noindent {\bf Definition 6.3.}
Let $M_1$ and $M_2$ be two metric spaces, with metrics $d_1$ and
$d_2$ respectively, and basepoints $x_1$ and $x_2$.
If $\epsilon$ is a positive real number,
then an {\sl $\epsilon$-approximation} between $(M_1, d_1, x_1)$ and
$(M_2, d_2, x_2)$ is a relation $R \subset M_1 \times M_2$ such that
\item{(i)} the projections $p_1 \colon R \rightarrow M_1$
and $p_2 \colon R \rightarrow M_2$ are both surjections,
\item{(ii)} if $xRy$ and $x'Ry'$, then
$\vert d_1(x,x') - d_2(y,y') \vert < \epsilon$, and
\item{(iii)} $x_1 R x_2$.

\noindent If $x$ is a point in a metric space $(M,d)$ and
$r$ is a positive real number, we denote the ball of radius $r$
about $x$ by $B_M(x,r)$. If we wish to emphasise the
metric on $M$, we may also write $B_{(M, d)}(x,r)$.
If $(M_\infty, d_\infty, x_\infty)$ and 
$\lbrace (M_i, d_i, x_i) \colon i \in {\Bbb N} \rbrace$
are metric spaces with basepoints, then
we say that the sequence $(M_i, d_i, x_i)$ {\sl converges} to
$(M_\infty, d_\infty, x_\infty)$ if,
for all $r > 0$, there is a sequence of positive real numbers
$\epsilon_i \rightarrow 0$ and $\epsilon_i$-approximations
between $B_{M_i}(x_i,r)$ and $B_{M_\infty}(x_\infty, r)$.
In this case, we write $(M_i, d_i, x_i) \rightarrow
(M_\infty, d_\infty, x_\infty)$.

The following example will be useful. Its proof (which is omitted)
is a straightforward application of Jacobi fields.

\noindent {\bf Lemma 6.4.} {\sl Let $M_i$ be a sequence
of simply-connected $n$-manifolds with
complete Riemannian metrics $g_i$ such that
$\kappa_{\rm \sup} (M_i, g_i) \rightarrow -1$
and $\kappa_{\rm \inf} (M_i, g_i) \rightarrow -1$.
Let $x_i$ be a basepoint in $M_i$, and
let $x_\infty$ be any point in hyperbolic $n$-space ${\Bbb H}^n$. 
Let $r$ be any positive real number. Then
for $i$ sufficiently large, there is a sequence
of real numbers $k_i > 1$ tending to $1$, and a sequence
of $k_i$-bi-Lipschitz homeomorphisms $h_i \colon B_{M_i}(x_i, r) \rightarrow
B_{{\Bbb H}^n}(x_\infty,r)$. In particular,
$(M_i, g_i, x_i)$ converges to ${\Bbb H}^n$ with
basepoint $x_\infty$.}

The following theorem of Gromov is a uniqueness result for
convergent sequences. A simple proof of this result can be found in [8].

\noindent {\bf Theorem 6.5.} [8] {\sl Let $(M_i, d_i, x_i)$ be
a sequence of complete metric spaces with basepoints, such that every
closed ball in $M_i$ is compact. Suppose that there are complete pointed metric
spaces $(M ,d, x)$ and $(M', d', x')$ such that
$$\eqalign {(M_i, d_i, x_i) &\rightarrow (M, d, x)\cr
(M_i, d_i, x_i) &\rightarrow (M', d', x').}
$$
Then $(M, d, x)$ and $(M', d', x')$ are isometric pointed metric spaces.}

Vital in our construction of a hyperbolic metric on ${\rm int}(M)$
will be the following theorem.

\noindent {\bf Theorem 6.6.} [8] {\sl Let $(M_i, d_i, x_i)$ be a sequence
of complete metric spaces (with basepoints)
in which bounded balls are compact. Then the following are equivalent.
\item{(1)} There is a subsequence $\lbrace
(M_j, d_j, x_j) \colon j \in J \subset
{\Bbb N} \rbrace$ converging to a complete metric space $(M, d, x)$.
\item{(2)} There is a subsequence $\lbrace (M_k, d_k, x_k)
\colon k \in K \subset
{\Bbb N} \rbrace$ such that for all $\epsilon > 0$ and $r > 0$, there is
a natural number $K(r, \epsilon)$ with the property that the
number of $\epsilon$-balls required to cover
$B_{M_k}(x_k, r)$ is less than $K(r, \epsilon)$.

}

Gromov also proved that, under certain circumstances, convergence
in the above sense implies bi-Lipschitz convergence.
The proof of this result [5] readily gives the following theorem.

\noindent {\bf Theorem 6.7.} {\sl Let $(M_i, g_i, x_i)$
be a sequence of Riemannian manifolds converging to some
space $(M_\infty, d_\infty, x_\infty)$. Let $r$ be a positive
real number. Suppose that
$\kappa_{\rm sup}(B_{M_i}(x_i,r), g_i)$ and 
$\kappa_{\rm inf}(B_{M_i}(x_i,r), g_i)$
are bounded above and below by constants independent of $i$.
Suppose also that the injectivity radius of $B_{M_i}(x_i, r)$
is bounded below by a constant independent of $i$. Then
for sufficiently large $i$, there is a sequence of real numbers
$k_i > 1$ tending to $1$, a sequence of positive real
numbers $\epsilon_i$ tending to zero and a sequence of 
$k_i$-bi-Lipschitz homeomorphisms
$h_i \colon B_{M_\infty}(x_\infty, r) \rightarrow U_i$,
where $B_{M_i}(x_i, r-\epsilon_i) \subset U_i \subset
B_{M_i}(x_i, r + \epsilon_i)$.}

The following lemma follows immediately from Theorem 6.6.

\noindent {\bf Lemma 6.8.} {\sl Let $(M_i, g_i, x_i)$ be a sequence
of $\delta_i$-pinched complete Riemannian $n$-manifolds
with $\delta_i \rightarrow 0$. Then some subsequence converges
to a pointed metric space $(M_\infty, d_\infty, x_\infty)$.}

\noindent {\sl Proof.} Let $(\tilde M_i, \tilde g_i)$ be the
universal cover of $(M_i, g_i)$ and let $\tilde x_i \in \tilde M_i$
be a lift of the basepoint $x_i$. Then Lemma 6.4
states that $(\tilde M_i, \tilde g_i, \tilde x_i)$ converges.
Hence, Theorem 6.6 states that we may pass to
a subsequence $(\tilde M_k, \tilde g_k, \tilde x_k)$ so that
for any $\epsilon > 0$
and $r>0$, there is a natural number $K(r, \epsilon)$ such that the
number of $\epsilon$-balls required to cover
$B_{\tilde M_k}(\tilde x_k, r)$ is less than $K(r, \epsilon)$.
Such a covering projects to a covering of $B_{M_k}(x_k, r)$
by $\epsilon$-balls. Hence, some subsequence
$(M_j, g_j, x_j)$ converges. $\square$

So, in our case where $M_i$ is a fixed manifold $M$,
some subsequence $(M, g_i, x_i)$ converges to a limit
$(M_\infty, d_\infty, x_\infty)$. In general, we lose
a great deal of information when passing from $(M, g_i)$
to $(M_\infty, d_\infty)$. In particular, $(M_\infty, d_\infty)$
need not be a hyperbolic 3-manifold. However, if
we pick the basepoints $x_i$ judiciously, then
$M_\infty$ will be hyperbolic. We shall pick $x_i$ in the
`thick' part of $(M, g_i)$. If $\epsilon$ is a
positive real number, then we denote the $\epsilon$-thick
part of $(M, g_i)$ by $(M, g_i)_{[\epsilon, \infty)}$ and
$\epsilon$-thin part of $(M, g_i)$ by $(M, g_i)_{(0, \epsilon]}$.

The Margulis lemma [3] describes the $\epsilon$-thin part of
hyperbolic manifolds for $\epsilon$ sufficiently small.
There is an extension of this result to negatively curved
Riemannian manifolds [2]. This implies that there is a positive real
number $\mu$ (called a Margulis constant) with the 
following property. If $M$ is
an orientable 3-manifold and $g$ is a $\delta$-pinched Riemannian
metric on $M$ with $\delta < 1$, then each
component $X$ of $(M,g)_{(0, \epsilon]}$ for $\epsilon \leq
\mu$ is diffeomorphic to one of the following possibilities.
\item{(i)} $X \cong D^2 \times S^1$. In this case, $X$ is known
as a `tube'. It is a neighbourhood of a closed geodesic
in $M$ with length less than $\epsilon$.
\item{(ii)} $X \cong S^1$. Then $X$ is a closed geodesic
with length precisely $\epsilon$. By perturbing our choice
of $\epsilon$ a little, we can ensure that this possibility
never arises.
\item{(iii)} $X \cong S^1 \times S^1 \times [0,\infty)$. Then
$X$ is a `cusp'.

In fact, a closer examination of the metric on the
$\mu$-thin part of $(M, g)$ readily yields the following
two results.

\noindent {\bf Proposition 6.9.} {\sl There is a
function $D \colon (0,\mu/2) \rightarrow {\Bbb R}_+$ with 
$D(\epsilon) \rightarrow \infty$ as $\epsilon \rightarrow 0$,
and which has the following property.
Pick real numbers $\delta$ and $\epsilon$
with $0 < \delta < 1$ and $0 < \epsilon < \min \lbrace 1, \mu/2 \rbrace$. 
If $M$ is a compact orientable 3-manifold
with a complete finite volume $\delta$-pinched Riemannian
metric $g$ on its interior, then the distance between
$(M,g)_{(0, 2\epsilon^2]}$ and $(M,g)_{[2\epsilon, \infty)}$
is at least $D(\epsilon)$.}

\noindent {\bf Proposition 6.10.} {\sl There is
a function $H \colon (0, \mu) \rightarrow {\Bbb R}_+$ with
$H(\epsilon) \rightarrow \infty$ as $\epsilon \rightarrow 0$
and which has the following property. Pick $\delta \in (0,1)$.
Let $M$ be a compact
orientable 3-manifold with a complete finite volume
$\delta$-pinched Riemannian metric $g$ on its interior.
If $X$ is a component of
$(M,g)_{(0, \mu]}$ which is a neighbourhood of a geodesic
of length at most $\epsilon$, then the length
of a meridian curve on $\partial X$ is at least
$H(\epsilon)$.}

We are now ready to prove Theorem 6.2.

\noindent {\sl Proof of Theorem 6.2.}
Since $M$ is almost hyperbolic, there is a sequence of positive
real numbers $\delta_i$ tending to zero,
and complete finite volume Riemannian metrics $g_i$ on ${\rm int}(M)$,
such that $(M,g_i)$ is $\delta_i$-pinched. 
We may assume that, for all $i$, $\delta_i < 1$. For each $i$,
pick a basepoint $x_i$ in $(M, g_i)_{[\mu, \infty)}$,
where $\mu$ is the constant mentioned above.

\vfill\eject

\noindent {\sl Claim 1.}  Some subsequence
$(M, g_j, x_j)$ converges to $(M_\infty, g_\infty, x_\infty)$,
which is a $3$-manifold $M_\infty$ with a complete hyperbolic
Riemannian metric $g_\infty$.

Some subsequence converges to a complete metric
space $(M_\infty, d_\infty, x_\infty)$ by Lemma 6.8.
Pass to this subsequence. Let $y$ be any point in $M_\infty$.
We wish to examine a neighbourhood of $y$.

Let $r$ be $d_{M_\infty}(x_\infty, y) + \mu$.
By definition, 
there is a sequence $\epsilon_i \rightarrow 0$
and $\epsilon_i$-approximations between $(B_{(M, g_i)}(x_i,r), g_i, x_i)$
and $(B_{M_\infty}(x_\infty, r), d_\infty, x_\infty)$.
From this, we get a sequence of $\epsilon_i$-approximations
between $(B_{(M, g_i)}(x_i,r), g_i, y_i)$ and
$(B_{M_\infty}(x_\infty, r), d_\infty, y)$ for some $y_i \in M$.
For any real number $z$ with $\epsilon_i < z < \mu$,
each $\epsilon_i$-approximation
restricts to an $\epsilon_i$-approximation between
$(U_i, g_i, y_i)$ and $(B_{M_\infty}(y, z), d_\infty, y)$,
where $B_{(M, g_i)}(y_i, z - \epsilon_i) \subset U_i
\subset B_{(M, g_i)}(y_i, z +\epsilon_i)$.
We can extend this to a $2\epsilon_i$-approximation between
$(B_{(M, g_i)}(y_i, z), g_i, y_i)$ and \hfill\break
$(B_{M_\infty}(y, z), d_\infty, y)$. So,
$(B_{(M, g_i)}(y_i, z), g_i, y_i) \rightarrow
(B_{M_\infty}(y,z), d_\infty, y)$.

If we insist that $\epsilon_i < \epsilon \leq \mu/2$,
then $x_i \in (M, g_i)_{[\mu, \infty)} \subset
(M, g_i)_{[2\epsilon, \infty)}$. By construction,
$y_i \in B_{(M, g_i)}(x_i, r)$. For $\epsilon$ sufficiently small, 
the distance between
$(M, g_i)_{(0, 2\epsilon^2]}$ and $(M, g_i)_{[2\epsilon, \infty)}$
is more than $r$, by Proposition 6.9. Hence,
$y_i \in (M, g_i)_{[2\epsilon^2,\infty)}$. Thus,
$B_{(M, g_i)}(y_i, \epsilon^2)$ is isometric to a ball of
radius $\epsilon^2$ in the universal cover $(\tilde M, \tilde g_i)$ of
$(M, g_i)$.
But $(\tilde M, \tilde g_i)$ is a complete simply-connected
Riemannian 3-manifold, and both
$\kappa_{\rm \sup} (\tilde M, \tilde g_i)$
and $\kappa_{\rm \inf} (\tilde M, \tilde g_i)$ tend to $-1$
as $i \rightarrow \infty$.
Thus, Lemma 6.4 states that $B_{(M,g_i)}(y_i, \epsilon^2)$ converges
to a ball of radius $\epsilon^2$ in ${\Bbb H}^3$, with
basepoint at the centre $p$ of the ball. Theorem 6.5
implies that this ball is isometric to $B_{M_\infty}(y, \epsilon^2)$,
via an isometry taking $p$ to $y$. Thus, $(M_\infty, d_\infty)$
is a 3-manifold with a complete hyperbolic Riemannian metric $g_\infty$. 
This proves the claim.

\noindent {\sl Claim 2.} Let $r$ be any positive real number.
For $i$ sufficiently large, there is a sequence of real numbers
$k_i > 1$ tending to $1$, a sequence of positive real
numbers $\epsilon_i$ tending to zero and a sequence of 
$k_i$-bi-Lipschitz homeomorphisms
$h_i \colon B_{M_\infty}(x_\infty, r) \rightarrow U_i$,
where $B_{(M, g_i)}(x_i, r-\epsilon_i) \subset U_i \subset
B_{(M, g_i)}(x_i, r + \epsilon_i)$.

We just need to check that the conditions of Theorem 6.7
are satisfied. The restrictions on 
$\kappa_{\rm sup}(B_{(M, g_i)}(x_i,r), g_i)$ and
$\kappa_{\rm inf}(B_{(M, g_i)}(x_i,r), g_i)$ hold automatically
since $(M, g_i)$ is $\delta_i$-pinched with $\delta_i \rightarrow 0$.
Since $x_i \in (M, g_i)_{[\mu, \infty)}$,
Proposition 6.9 implies that there is some $\epsilon \leq \mu/2$
such that $(M, g_i)_{[2\epsilon^2, \infty)} \supset
B_{(M, g_i)}(x_i, r)$ for all $i$ sufficiently large.
Thus Theorem 6.7 proves the claim.

\noindent {\sl Claim 3.} $(M_\infty, g_\infty)$ has finite volume.

In the proof of Proposition 3.3, we established that
$$\vert M \vert \pi/2 \geq (-\kappa_{\rm sup}(M, g_i))^{3/2}
{\rm Vol}(M, g_i).$$
Since $\kappa_{\rm sup}(M, g_i) \rightarrow -1$, we deduce that
the sequence ${\rm Vol}(M, g_i)$ is bounded above.
Now,
$${\rm Vol}(M_\infty, g_\infty) = \lim_{r \rightarrow \infty}
{\rm Vol}(B_{M_\infty}(x_\infty, r), g_\infty),$$
and, using the notation in Claim 2,
$$\eqalign{
{\rm Vol}(B_{M_\infty}(x_\infty, r), g_\infty) &\leq
(k_i)^3 {\rm Vol}(B_{(M, g_i)}(x_i, r + \epsilon_i), g_i)\cr
&\leq (k_i)^3 {\rm Vol}(M, g_i).}$$
Thus, ${\rm Vol}(M_\infty, g_\infty)$ is finite.

This implies that $M_\infty = {\rm int}(\overline M_\infty)$
for some compact orientable 3-manifold $\overline M_\infty$,
with $\partial \overline M_\infty$ a (possibly empty) union
of tori.

\noindent {\sl Claim 4.} Let $\epsilon$ be a positive real number
less than $\mu$ such that 
$(M_\infty, g_\infty)_{(0, \epsilon]}$ is either
empty or consists
only of horoball neighbourhoods of cusps. Then, for $i$ sufficiently large,
there is a sequence of real numbers $k'_i > 1$ tending
to $1$, and $k'_i$-bi-Lipschitz homeomorphisms
$h'_i \colon (M_\infty, g_\infty)_{[\epsilon, \infty)} \rightarrow
(M, g_i)_{[\epsilon, \infty)}$.

We pick $r > 0$ so that $B_{M_\infty}(x_\infty, r) \supset
(M_\infty, g_\infty)_{[\epsilon, \infty)}$.
Let $T$ be the boundary of 
$(M_\infty, g_\infty)_{[\epsilon, \infty)}$ which
is a collection of tori.
Using the notation in Claim 2, the homeomorphism
$h_i \colon B_{M_\infty}(x_\infty, r) \rightarrow U_i$
is almost an isometry for $i$ large. Therefore
for large $i$, there is a sequence of positive real numbers
$\gamma_i$ tending to zero, such that $h_i(T)$
separates $(M, g_i)_{(0, \epsilon - \gamma_i]}$
from $(M, g_i)_{[\epsilon + \gamma_i, \infty)}$.
But, using the Margulis Lemma for negatively curved
3-manifolds, $(M, g_i)_{[\epsilon - \gamma_i, \epsilon+\gamma_i]}$
is homeomorphic to a collection of copies of $T^2 \times I$.
Hence, $h_i(T)$ is isotopic to the boundary of $(M, g_i)_{[\epsilon, \infty)}$,
since any torus in $T^2 \times I$ which separates the
boundary components is isotopic to either boundary component.
We may therefore modify $h_i$ to $h'_i \colon
(M_\infty, g_\infty)_{[\epsilon, \infty)} \rightarrow
(M, g_i)_{[\epsilon, \infty)}$, ensuring that the $h'_i$
are bi-Lipschtiz homeomorphisms as claimed.

\noindent {\sl Claim 5.}  Let $\epsilon$ be a positive real number
less than $\mu$ such that 
$(M_\infty, g_\infty)_{(0, \epsilon]}$ is
either empty or consists
only of horoball neighbourhoods of cusps.
Then the length of the core geodesic of
each tube component of $(M, g_i)_{(0, \epsilon]}$ tends to
zero, as $i \rightarrow \infty$.

If not, we may find a positive real number $\alpha \leq \epsilon$ and 
a subsequence in which $(M, g_i)_{(0, \epsilon]}$
contains a geodesic of length at least $\alpha$.
Applying Claim 4, we find that, for $i$ sufficiently large,
there is a sequence of real numbers $k''_i > 1$ tending
to $1$, and $k''_i$-bi-Lipschitz homeomorphisms
$h''_i \colon (M_\infty, g_\infty)_{[\alpha, \infty)} \rightarrow
(M, g_i)_{[\alpha, \infty)}$. In particular, there is
a geodesic in $(M_\infty, g_\infty)$ of length at
most $\alpha k''_i$. But $\alpha \leq \epsilon$,
which is less than the length of the shortest geodesic
in $M_\infty$. This is a contradiction, which proves the claim.

Fix $\epsilon \leq \mu$ such that 
$(M_\infty, g_\infty)_{(0, \epsilon]}$ is
either empty or consists
only of horoball neighbourhoods of cusps.
Now, $(M, g_i)_{(0, \epsilon]}$ is a (possibly empty) collection
of tubes and a (possibly empty) collection of cusps.
If, for infinitely many $i$, $(M, g_i)_{(0, \epsilon]}$ contains no tubes, then
Claim 4 implies that $M$ is homeomorphic to $\overline M_\infty$.
By Claims 1 and 3, $M_\infty$ has a complete finite volume
hyperbolic structure, which proves the theorem in this case.

Consider now the case where $(M, g_i)_{(0, \epsilon]}$ is, for
infinitely many $i$, a collection of
tubes $X^1_i, \dots X^{n(i)}_i$ and possibly also
some cusps. Claim 4 implies that, for each $i$ sufficiently large,
the meridian slope on $X^j_i$ corresponds to
a slope $s^j_i$ on $\partial \overline M_\infty$, and that $M$
is homeomorphic to $\overline M_\infty(s^1_i, \dots, s^{n(i)}_i)$.
By passing to a subsequence, we may ensure that
$n(i)$ is some fixed integer $n$, and that, for each $j$,
the slopes $s^j_i$ all lie a fixed torus $T_j$.
Claim 5 states that the core geodesic of $X^j_i$
tends to zero as $i \rightarrow \infty$. Hence, Proposition 6.10
states that the length of the meridian slope on
$X^j_i$ tends to infinity. 
The length of $s^j_i$ on $(M_\infty, g_\infty)_{[\epsilon, \infty)}$
differs from length of the corresponding meridian
slope on $X^j_i$ by a factor of at most
$k'_i$, which converges to $1$, as in Claim 4.
Thus, $l(s^j_i) \rightarrow \infty$
as $i \rightarrow \infty$. Therefore, by passing to
a subsequence, we may assume that the slopes $s^j_i \not=
s^j_k$ if $i \not= k$. Lemma 4.2 gives us a contradiction.
$\square$

\noindent {\bf Remark.} The proof of Theorem 6.2 actually gives
something a little stronger. It shows that if $(M, g_i)$
is a sequence of $\delta_i$-pinched Riemannian manifolds
with $\delta_i \rightarrow 0$, then ${\rm int}(M)$ has a complete
finite volume hyperbolic metric $h$, and (for all $i$
sufficiently large) there is a sequence
of real number $k_i > 1$, tending to 1, and a sequence
of $k_i$-bi-Lipschitz homeomorphisms between $(M, g_i)$
and $(M, h)$.

\vskip 18pt
\centerline{REFERENCES}
\vskip 6pt

\item{1.} C. ADAMS, {\sl The noncompact hyperbolic 3-manifold of minimal
volume,} Proc. A.M.S. {\bf 100} (1987) 601-606.
\item{2.} W. BALLMAN, M. GROMOV and V. SCHROEDER, {\sl
Manifolds of Non-positive Curvature}, Birkh\"auser (1985)
\item{3.} R. BENEDETTI and C. PETRONIO, {\sl Lectures on Hyperbolic Geometry},
\break
Springer-Verlag (1992).
\item{4.} S. BLEILER and C. HODGSON, {\sl Spherical space forms and Dehn
filling}, Topology {\bf 35} (1996) 809-833.
\item{5.} R. GREEN and H. WU, {\sl Lipschitz convergence of Riemannian
manifolds}, Pacific J. Math. {\bf 131} (1988) 119-141.
\item{6.} M. GROMOV, {\sl Structures m\'etriques pour les vari\'et\'es
Riemanniennes}, Cedic-Fernand-Nathan (1981)
\item{7.} H. HILDEN, M. LOZANO and J. MONTESINOS, {\sl On knots
that are universal}, Topology {\bf 24} (1985) 499-504.
\item{8.} C. HODGSON, {\sl Notes on the Orbifold Theorem}.
\item{9.} R. KIRBY, {\sl Problems in Low-Dimensional Topology}, U.C. Berkeley
(1995).
\item{10.} M. LACKENBY, {\sl Surfaces, surgery and unknotting operations},
Math. Ann. {\bf 308} (1997) 615-632.
\item{11.} M. LACKENBY, {\sl Dehn surgery on knots in 3-manifolds}, J.
Amer. Math. Soc. {\bf 10} (1997) 835-864.
\item{12.} J. MORGAN and H. BASS, {\sl The Smith Conjecture}, Academic Press
(1984).
\item{13.} R. MYERS, {\sl Open book decompositions of 3-manifolds},
Proc. Symp. Pure Math., AMS {\bf 32} (1977) 3-6.
\item{14.} P. PETERSEN, Private communication.
\item{15.} D. ROLFSEN, {\sl Knots and Links}, Publish or Perish (1976).
\vfill\eject
\item{16.} P. SCOTT, {\sl The geometries of 3-manifolds}, Bull. London
Math. Soc. {\bf 15} (1983) 401-487.
\item{17.} W. THURSTON, {\sl The geometry and topology of 3-manifolds},
Princeton University (1979).
\item{18.} Q. ZHOU, {\sl On the topological type of 3-dimensional
negatively curved manifolds}, Adv. in Math. (China) {\bf 22} (1993)
270-281.

\vskip 12pt
\+ Mathematics Department\cr
\+ University of California \cr
\+ Santa Barbara, CA 93106 \cr
\+ U.S.A. \cr
\vskip 12pt
\+ Department of Pure Mathematics and Mathematical Statistics \cr
\+ University of Cambridge \cr
\+ 16 Mill Lane \cr
\+ Cambridge CB2 1SB \cr
\+ England \cr

\end